\documentclass[11pt]{amsart}
\allowdisplaybreaks[1]

\usepackage[usenames]{color}

\usepackage{a4wide}
\usepackage{amssymb}
\usepackage{amsmath} 
\usepackage{amscd} 
                
\theoremstyle{plain}
\newtheorem{theorem}{Theorem} 
\newtheorem{lemma}{Lemma}
\newtheorem{prop}{Proposition} 
\newtheorem{coro}{Corollary}
\newtheorem{fact}{Fact} 

\theoremstyle{definition}
\newtheorem{definition}{Definition}
\newtheorem{remark}{Remark}

\newcommand{\ts}{\hspace{0.5pt}}

\newcommand{\CC}{\mathbb{C}\ts}
\newcommand{\RR}{\mathbb{R}\ts}
\newcommand{\QQ}{\mathbb{Q}\ts}
\newcommand{\ZZ}{\mathbb{Z}}
\newcommand{\NN}{\mathbb{N}}

\newcommand{\bA}{\mathbb A}
\newcommand{\TT}{\mathbb T}
\newcommand{\XX}{\mathbb X}
\newcommand{\EE}{\mathbb E}
\newcommand{\bS}{\mathbb S}
\newcommand{\KK}{\mathbb K\ts}
\newcommand{\LL}{\mathcal L}

\newcommand{\ba}{\beta^{}_{\bA}}
\newcommand{\bt}{\beta^{}_{\TT}}
\newcommand{\be}{\beta^{}_{\EE}}

\newcommand{\gL}{\varLambda}

\newcommand{\gG}{\varGamma}
\newcommand{\oplam}{\mbox{\Large $\curlywedge$}}

\newcommand{\Ghat}{\widehat{G}}
\newcommand{\gammahat}{\widehat{\gamma}}

\newcommand{\CalD}{\mathcal{D}}
\newcommand{\CalDq}{\mathcal{D}^{\equiv}}

\newcommand{\dd}{\,\mathrm{d}}
\newcommand{\symdiff}{\bigtriangleup}
\newcommand{\supp}{\mathrm{supp}}
\newcommand{\dens}{\mathrm{dens}}

\begin{document}

\title[Characterization of model sets]
{Characterization of model sets  \\[1mm] 
by dynamical systems}

\author{Michael Baake}
\address{Fakult\"at f\"ur Mathematik,
Universit\"at Bielefeld, Postfach 100131,
33501 Bielefeld, Germany}
\email{mbaake@math.uni-bielefeld.de}
\urladdr{http://www.math.uni-bielefeld.de/baake/}

\author{Daniel Lenz}
\address{Fakult\"at f\"ur Mathematik, TU Chemnitz,
09107 Chemnitz, Germany}
\email{dlenz@mathematik.tu-chemnitz.de}
\urladdr{http://www.tu-chemnitz.de/mathematik/analysis/dlenz}

\author{Robert V.\ Moody}
\address{Department of Mathematics and Statistics, 
University of Victoria, \newline
\hspace*{12pt}Victoria, British Columbia V8W 3P4, Canada}
\email{rmoody@uvic.ca}
\urladdr{http://www.math.ualberta.ca/\~{}rvmoody/rvm/}

\begin{abstract} 
  It is shown how regular model sets can be characterized in terms of
  regularity properties of their associated dynamical systems. The
  proof proceeds in two steps. First, we characterize regular model
  sets in terms of a certain map $\beta$ and then relate the
  properties of $\beta$ to ones of the underlying dynamical system. As
  a by-product, we can show that regular model sets are, in a suitable
  sense, as close to periodic sets as possible among repetitive
  aperiodic sets.
\end{abstract}

\maketitle

\section{Introduction}
 
Delone sets provide an important model class for the description of
aperiodic order. In particular, they can be viewed as a mathematical
abstraction of the set of atomic positions of a physical quasicrystal
(at zero temperature, or at a given instant of time). Many of the
rather intriguing spectral properties of quasicrystals can be
formulated, in a simplified manner, on the basis of Delone sets.  The
latter contain the important class of model sets (see below for
definitions), which is our main topic in this paper.

Since the discovery of quasicrystals \cite{Shechtman}, model sets have
been a particular focus of attention because they are, except under
extreme conditions, pure point diffractive
\cite{Hof,Moody,Martin2}. This property remains true also under
certain equivariant perturbations, which turn them into deformed model
sets \cite{Hof,BD,BL2}, and extend the applicability of these sets
considerably \cite{SW}.

Model sets are discrete point sets that arise by (partial) projection
of a lattice from some ``higher dimensional" or ``super" space. To
avoid misunderstandings, and to accommodate situations where the
concept of dimension is not available, we shall call this super space
the \emph{embedding space} below. Model sets have been found useful
in numerous studies both by experimentalists modelling quasicrystals
and by mathematicians studying aperiodic order and diffraction. One
principal difficulty has been to find good characterizations of them.
In particular, what are the intrinsic properties of a point set that
permits its description as a projection from (parts of) some higher
dimensional lattice?

Another major ingredient in the study of aperiodic point sets (and
tilings) has been the use of dynamical systems.  Given a (suitably
discrete) point set $\gL \subset \RR^d$, for example, one associates
with it a space which is the closure of its $\RR^d$-translation orbit,
the closure taken in a topology that compares point sets for more or
less exact match in local regions around the origin. This is called
the dynamical hull, or \emph{local hull} in this paper (since we shall
meet other hulls that are dynamical systems as well).  The major
objective of this paper is to characterize model sets in terms of the
properties of their local hulls.

As the theory of model sets and related mathematics has developed, it
has become clear that the properties of the ambient space that are
required are sufficiently weak that the group $\RR^d$ may be replaced
by any $\sigma$-compact locally compact Abelian (LCA) group $G$,
without increasing the complexity of the proofs. In fact, this
additional generality is to some extent necessary to understand model
sets, as we shall see. In this paper, we take this more general setting.
 
The main theorem of the paper is:
 
\begin{theorem}\label{main1} 
   Let $G$ be a $\sigma$-compact LCA group and\/ $(\XX,G)$ a point set
   dynamical system on $G$.  Then, for $(\XX,G)$ to be the dynamical
   system associated to a repetitive regular model set it is necessary
   and sufficient for the following four conditions to be satisfied.
\begin{itemize}
\item[(1)] All elements of\/ $\XX$ are Meyer sets;
\item[(2)] $(\XX,G)$ is strictly ergodic, i.e., uniquely ergodic and
           minimal;
\item[(3)] $(\XX,G)$ has pure point dynamical spectrum with continuous
           eigenfunctions;
\item[(4)] The eigenfunctions of $(\XX,G)$ separate almost all points
           of\/ $\XX$ $($i.e., the set\newline
           $\{\gG \in \XX : \mbox{there exists $\gG'\neq \gG$ with
           $f(\gG) = f(\gG')$ for all eigenfunctions $f$}\}$ \newline
           has measure $0 )$.
\end{itemize}
\end{theorem}

The necessity of the conditions is already known
\cite{Martin1,Martin2}, so our main task is to deal with the converse
-- the four listed properties characterize repetitive regular model
sets --, although we end up proving the necessity again in the
process.

The proof is broken into three main parts. The first part is to use the
properties (3) and (4) to identify elements of $\XX$ that cannot be
separated by the continuous eigenfunctions. This results in a new
dynamical system $(\EE,G)$, where $\EE$ is actually a compact Abelian group,
and a surjective $G$-mapping of $\XX$ onto $\EE$. Although this new
group need not be a torus, it is nonetheless useful to simply call
such a map a torus map or \emph{torus parametrization}, in analogy to
\cite{BHP}.

The second part is to show that $\EE$ can be identified with
another dynamical hull $\bA$ of $\gL$, this time determined not by the
local topology, but rather by a topology called the
\emph{autocorrelation topology}. This topology compares point sets
globally for statistical match.

The third step, which actually appears first in the paper, is to show
that a torus mapping of $\XX$ onto $\bA$ assures that we are in the
situation of model sets -- we can explicitly construct the embedding
space, the lattice, and the mechanism which controls the projection
down into the ambient space. This is really the heart of the
matter. Given a Meyer set $\gL$, we have its two hulls $\XX(\gL)$
and $\bA(\gL)$. These are quite natural objects.  The mapping $\beta$
between them, when it exists, is the most natural one possible. It is
really nothing but looking at the same elements of $\XX(\gL)$, but in
another topology. The assumption of the existence of the map is the
same as saying that this change of topology is continuous, which in
turn is the same as requiring that the local and global topologies are
consistent with each other. It is this consistency that effectively
characterizes the cut and project formalism.

The existence of windows for realizing the elements of $\XX$ as model
sets (or inter model sets) emerges as we require more out of the
mapping $\beta\,$: first that it is one-to-one somewhere, and finally that it
is one-to-one almost everywhere.  If we go so far as to assume that it is
one-to-one everywhere, we collapse into the crystallographic case
(Theorem~\ref{main5}). Thus condition (4) of Theorem~\ref{main1} seems
to contain the essence of aperiodicity (at least in the context of
Meyer sets).  This gives another instance for the intuition that
regular model sets are a very natural generalization of
crystallographic (i.e., fully periodic) point sets, and that aperiodic
model sets are, in this sense, as close to periodic sets as possible
among (repetitive) aperiodic Meyer sets.

Section \ref{definitions} introduces the basic definitions and
concepts used throughout the paper. In particular, in
Paragraphs~\ref{basic}, \ref{local hull} and \ref{ac hull}, we
establish the basic notions about the point sets and dynamical hulls
that appear in the paper. Paragraph \ref{model} deals with
cut and project schemes and model sets. Paragraph~\ref{torus}
introduces the notion of a torus parametrization. While the material of
that section is essentially known, the point of view taken there is of
fundamental importance for our considerations.

Beyond the main theorem, there are a number of intermediate results
that are interesting in their own right and are part of the overall
proof.  Section~\ref{results} serves the purpose of detailing these
results and indicating the logical flow of the paper. The paper proper
then begins with the consequences of a torus parametrization $\beta \!:\,
\XX(\gL) \rightarrow \bA(\gL)$, gradually refining what can be learned
from it as further conditions are added.

Model sets, as one sees them in the literature, come with varying
definitions and side conditions, depending on the requirements of the
moment. However, our results require quite precise notions of what
constitutes a cut and project scheme, which windows are permitted, and
how they relate to the cut and project scheme. Much, but not all, of
this appears in the work of Schlottmann cited above. To make things
clear, particularly the important ideas of irredundancy, which is not
standard, and inter model sets, which are new~\cite{LM}, we have
reworked this material and included it in the paper. Our attitude is
that the main purpose of the paper is to prove the sufficiency
direction of the Theorem~\ref{main1}, whence we have written the paper
so that it moves in that direction from the very beginning.
By the time that we have proved sufficiency, we actually know enough
to prove necessity rather easily.

The paper has been delayed in reaching its final form by various
circumstances around the lives of its authors. Nonetheless, its
results have been announced in several places
\cite{MoodyLisbon,MoodyOberwolfach}. Meanwhile, based on this
paper, an extension of part of this theory to multi-colour sets has
been worked out \cite{LM}, and this has been effectively used in
establishing the equivalence of pure pointedness and model sets for
substitution tilings and point sets \cite{Lee}, a result that, for the
case of unimodular Pisot substitutions in one dimension, has recently
also been discussed in a slightly different context \cite{K}.

\section{Basic definitions and hulls} \label{definitions}

This paper is a study of the relationship between various concepts
in the regime of aperiodic order, formulated in terms of point sets in
locally compact Abelian (LCA) groups. Let us first introduce the
concepts.

\subsection{Aperiodic order and diffraction theory:\ the general setting}
\label{basic} 

Let $G$ be a locally compact Abelian group, with Haar measure
$\theta_G$ (normalized as $\theta_G (G)=1$ if $G$ is compact).  We
assume that $G$ is $\sigma$-compact (also called countable at
infinity).  This is equivalent to the existence of an \emph{averaging
sequence} $\mathcal{A} = \{ A_n : n\in\NN\}$ of open, relatively
compact sets $A_n \subset G$ with $\overline{A_n}\subset A_{n+1}$ and
$G=\bigcup_{n\ge 1} A_n$.  In fact, the averaging sequence can be
chosen to be a \emph{van Hove sequence}, see \cite{Martin2} for
details. This means that, for every compact set $K\subset G$,
\[
   \lim_{n\to\infty} \frac{\theta^{}_G \big( (( K + A^{}_n )\setminus
  A_n^\circ ) \cup (( - K +\overline{G\setminus A^{}_n})\cap A_n )
  \big)}{\theta^{}_G (A^{}_n)} \; = \; 0,
\]
where the bar (circle) denotes the closure (interior) of a set. In
effect, this rather technical looking condition states that for each
compact subset $K$ of $G$, the $K$-boundary of the averaging sequence
becomes negligible (in the sense of measure) to the sequence itself as
$n\to\infty$. Note that general van Hove sequences need not be nested.

\smallskip
A subset $\gL$ of $G$ is called \emph{$U$-uniformly discrete} if, for
the open neighbourhood $U$ of $0$ in $G$ and for all $x\in \gL$,
$(x+U) \,\cap \, \gL = \{x\}$. We say that $\gL$ is \emph{uniformly
discrete} if a neighbourhood $U$ exists for which $\gL$ is
$U$-uniformly discrete.  By $\sigma$-compactness of $G$, every
uniformly discrete set in $G$ is countable.  The set of all uniformly
discrete subsets of $G$ is denoted by $\CalD =\CalD (G)$ and the set
of $U$-uniformly discrete subsets by $\CalD_U$.

Uniformly discrete sets can have various further regularity
properties.  A uniformly discrete subset $\gL$ of $G$ is called
\emph{Delone} if it is also relatively dense, i.e., if there exists
a compact set $K$ in $G$ with $G = \gL + K$.

Now, let $\gL$ be an arbitrary uniformly discrete set. Then, $\gL$ is
of \emph{finite local complexity} (FLC) if the set of $K$-clusters,
\[   
    \{ (-x + \gL)\cap K : x\in \gL\}  \, , 
\]
is finite for every compact $K\subset G$. This is equivalent to 
$\gL - \gL$ being discrete and closed \cite{Martin2}.  
If $\gL \subset G$ is a Delone set and there exists a finite set
$F\subset G $ with $\gL -\gL\subset \gL + F$, then $\gL$ is called a
\emph{Meyer set}.  Evidently, $\gL - \gL$ is uniformly discrete
whenever $\gL$ is Meyer.

A point set $\gL\subset G$ of finite local complexity is called
\emph{repetitive} if for every compact $K$ in $G$ the set of
repetitions of $\gL \cap K$
\[
   \{t\in G : (-t + \gL) \cap K  = \gL \cap K\}
\] 
is relatively dense in $G$. $\gL$ is said to have \emph{uniform patch
  frequencies} (some people say \emph{uniform cluster frequencies})
if, for each finite subset $P$ of $G$ and for all $a\in G$,
 \[
   \frac{ \mathrm{card} \{ t\in G : t+P \subset \gL \cap
   (a + B_n)\}} {\theta_G (B_n) }
\]
converges uniformly in $a$ along \emph{every} van Hove sequence
$\{B_n : n\in\NN\}$, see \cite{Martin2} for details.

\smallskip
The diffraction pattern of a solid modelled by $\gL$ can be described
as follows \cite{Cowley,Hof}. For $x\in G$, let $\delta_x$ denote the
normalized point (or Dirac) measure at $x\in G$. If the limit (taken
in the vague topology)
\[
   \gamma \; := \; \lim_{n\to \infty} \frac{1}{\theta^{}_{G}(A^{}_{n})}
             \sum_{x,y\in \gL\cap A_n} \delta_{x-y}
\]
exists, it is called the \emph{autocorrelation measure} of $\gL$
relative to the averaging sequence $\mathcal{A}$. If $\gL$ has uniform
patch frequencies, the limit exists and does not depend on the choice
of $\mathcal A$ (as long as it is van Hove).  The autocorrelation
measure is positive definite and hence transformable, i.e., we can
take its Fourier transform $\widehat{\gamma}$. This is a positive
measure on the dual group $\widehat{G}$, called the \emph{diffraction
  measure}. For $G=\RR^n$, it describes the outcome of a diffraction
experiment, compare \cite{Cowley} for details.

\subsection{The local hull}\label{local hull}

In this and the next paragraph, we introduce two topologies on the
set $\CalD$ of all uniformly discrete subsets of $G$. The interplay of
these two topologies is a main feature of the paper.

The so-called \emph{local topology} (LT) on $\CalD$ is defined via
the uniform structure given by the entourages
\[
   U^{}_\mathrm{LT} (K,V) \; := \; \{(\gG,\gG')\in\CalD\times\CalD : 
   (v + \gG)\cap K = \gG'\cap K \;\mbox{for some $v\in V$}\}
\]
for $K\subset G$ compact and $V$ a neighbourhood of $0$ in $G$.  Thus,
two uniformly discrete sets are close if they agree on a ``large''
compact set up to a ``small'' (global) translation. For definitions,
terminology and basic theorems on uniformities, see \cite{Bou,Q}.

As is immediate from the definition of the local topology, the
canonical action of $G$ on $\CalD$ given by 
\[ 
   G\times \CalD\longrightarrow \CalD, \quad 
   (t,\gL)\mapsto -t + \gL, 
\]
is continuous. In particular, if $\XX\subset \CalD $ is compact in the
local topology and invariant under this action, then $(\XX,G)$ is a
topological dynamical system. Such a dynamical system will be called a
\emph{point set dynamical system}.

The hull of an element $\gL\in \CalD$ in the local topology (i.e., the
closure of the orbit $ G + \gL = \{ x+\gL : x\in G\}$) is denoted by
$\XX (\gL)$. 
\begin{fact}\label{LT-compact} \cite{Martin2}
  If $\gL$ is a Delone set, the hull\/ $\XX(\gL)$ is LT-compact if and
  only if $\gL$ is of finite local complexity, i.e., if and only if\/
  $\gL-\gL$ is discrete and closed. \qed
\end{fact} 
In this case, $\XX (\gL)$ gives rise to a point set dynamical system
$(\XX (\gL),G)$.  This dynamical system is a basic object in the study
of the long-range order of discrete point sets because of its ability
to reflect important geometric properties in the language of
dynamical systems.

\begin{fact} \cite{Martin2}
   Let $\gL$ be a Delone set of finite local complexity.  Then, the
   dynamical system $(\XX(\gL),G)$ is uniquely ergodic\/ $($i.e.,
   there exists precisely one $G$-invariant probability measure on
   $\XX(\gL)$\/$)$ if and only if $\gL$ has uniform patch frequencies.
   \qed
\end{fact}

Two Delone sets $\gL,\gL'$ are \emph{locally indistinguishable} (LI) if
each cluster of $\gL$ (i.e., each set of the form $\gL \cap K$ with
$K\subset G$ compact) is a translate $-x + (\gL' \cap (x + K))$ of a
cluster of $\gL'$ and vice versa.
 This
equivalence relation defines the so-called LI classes, and one has
\begin{fact}
  Let $G$ be an LCA group and $\gL\subset G$ a Delone set of finite
  local complexity. Then, the following properties are equivalent.
\begin{itemize}
\item[{(1)}] The set $\gL$ is repetitive.
\item[{(2)}] The hull\/ $\XX (\gL)$ is the LI class of
$\gL$.
\item[{(3)}] The dynamical system $(\XX (\gL),G)$ is
minimal.
\end{itemize}
\end{fact}
\begin{proof} This is a variant of Gottschalk's theorem, see
\cite{Martin2} for details.
\end{proof}

The definition of closeness in the local topology has a special
consequence for translates of Meyer sets. 

\begin{fact} \label{snapin} 
  Let $\gL$ be a Meyer set.  Then, for all suitably small
  neighbourhoods $V$ of\/ $0$ in $G$ and all compact $C\subset G$ 
  with $C \cap \gL\neq \varnothing$, the equality $ \gL\cap C = (-x +
  \gL)\cap C $ holds whenever $(- x + \gL,\gL)\in U^{}_\mathrm{LT}
  (C,V)$ for $x\in \gL -\gL$.
\end{fact}
\begin{proof} As $\gL$ is Meyer, one has $\gL - \gL \subset \gL + F$
  with $F$ a finite set. Clearly, also $(\gL -\gL) + (\gL - \gL)
  \subset \gL + F'$, with $F'$ still finite, so that uniform
  discreteness persists to $\gL -\gL + (\gL - \gL)$. Thus, there
  exists an open neighbourhood $V$ of $0$ in $G$ so small that $V \cap
  ((\gL - \gL) + (\gL - \gL))= \{0\}$. Now, $(- x + \gL,\gL)\in
  U^{}_\mathrm{LT} (C,V)$ for $x\in G$ implies
\[ 
   (v - x + \gL) \cap C = \gL \cap C 
\] 
for some $v\in V$.  Now, if $x\in \gL-\gL$, then $\gL \cap C\neq
\varnothing $ yields $v\in V \cap ((\gL- \gL) + (\gL - \gL)) = \{0\}$
and the fact is proved.
\end{proof}

\smallskip
Let $(\XX,G)$ be a point set dynamical system which is uniquely
ergodic.  In this case, there is a canonical Hilbert space associated
to $(\XX,G)$, the space $L^2(\XX,\mu)$ of square integrable functions
on $\XX$ (with respect to the unique $G$-invariant probability measure
$\mu$).  The action of $G$ on $\XX$ gives rise to a unitary
representation $T$ of $G$ on this space via
\begin{equation*}
  T_t \! : \;  L^2(\XX,\mu)\; \longrightarrow \;L^2(\XX,\mu)  \, , 
  \quad (T_t f) (\gL) \; := \; f(-t + \gL) \, ,
\end{equation*}
for $f\in L^2(\XX,\mu) $ and $t\in G$.  An $f\in L^2(\XX,\mu)$ is
called an \emph{eigenfunction} of $T$ with \emph{eigenvalue}
$\hat{s}\in \Ghat$ (the dual group) if $T_t f = (\hat{s}, t) f$ for
every $t\in G$, where $(\hat{s}, .)$ denotes the character defined by
$\hat{s}$.  An eigenfunction (to $\hat{s}$, say) is called
\emph{continuous} if it has a continuous representative $f$ with $f(-
t + \gL ) = (\hat{s},t)\ts f(\gL)$, for all $\gL\in \XX$ and $t\in G$.
The representation $T$ is said to have \emph{pure point spectrum} if
the set of eigenfunctions is total in $L^2 (\XX,\mu)$.  One then also
says that the dynamical system $(\XX,G)$ has \emph{pure point
  dynamical spectrum}.

\subsection{The autocorrelation hull}\label{ac hull}

The \emph{upper density} of a point set $\gL \subset G$ is defined by
\[ 
    \overline{\dens} (\gL) \; := \; \limsup_{n\to \infty} 
    \frac{\mathrm{card} ( \gL \, \cap A^{}_n) }{\theta^{}_G(A^{}_n)}
\]
with respect to the averaging van Hove sequence $\mathcal A$ chosen
before. The lower density, $\underline{\dens} (\gL)$, is defined
analogously. If $\overline{\dens} (\gL)=\underline{\dens} (\gL)$, this
is called the \emph{density} of $\gL$, denoted by $\mathrm{dens}
(\gL)$. We shall usually suppress the explicit reference to
$\mathcal{A}$.

The \emph{mixed autocorrelation topology} (mACT) on $\CalD$ is defined via
the uniform structure given by the entourages
\[
   U^{}_\mathrm{mACT} (V,\varepsilon) \; := \; \{(\gG,\gG') \in 
         \CalD \times\CalD :  d(v + \gG,\gG')\leq \varepsilon  
         \;\mbox{for some $v\in V$}\},
\]
for every neighbourhood $V$ of $0$ in $G$ and every $\varepsilon >0$,
where the pseudo-metric $d$ on $\CalD$ is defined by the upper density
of the symmetric difference of sets:
\begin{equation}\label{metric}
   d(\gG,\gG') \; := \;
   \overline{\dens}\,(\gG \, \symdiff \, \gG') \, .
\end{equation} 

Note that the triangle inequality follows from the fact that
$\gG\symdiff\gG' \subset (\gG\symdiff\gG'') \cup (\gG'\symdiff\gG'')$,
for arbitrary point sets $\gG''\subset \mathcal{D}$.  With this
definition, $d$ is $G$-invariant (i.e., $d(t + \gG, t + \gG') = d(\gG,
\gG')$ for all $t\in G$ and all $\gG,\gG'\in \CalD$), because
$\mathcal{A}$ has the van Hove property.  We call mACT the
\emph{mixed} autocorrelation topology because it mixes the ordinary
topology of $G$ with the topology introduced by the pseudo-metric $d$.
The topology induced by $d$ itself, in turn, ultimately arises from
the autocorrelation (see below) and we thus call it the
\emph{autocorrelation topology}.

Note that $d$ contains information on statistical coincidence of the
global structure. Thus, two sets are close in the mixed
autocorrelation topology if their global structures are statistically
close up to a small translation.

\smallskip 

It is obvious that, for a general LCA group $G$, $d$ does
  not define a metric on $\mathcal{D}$. However, it does permit the
  construction of a completion where $d$ becomes a metric.  To see
this, fix an open neighbourhood $U$ of $0$ in $G$ and consider the
restriction of the pseudo-metric $d$ (still denoted by $d$) to
$\CalD_U$. Introduce an equivalence relation $\equiv$ on $\CalD_U$ by
setting $\gG \equiv \gG'$ if and only if $d(\gG,\gG')=0$, with $d$ as
defined in (\ref{metric}). The quotient of $\CalD_U$ by this
equivalence relation is denoted by $\CalDq_U$. By construction, the
pseudo-metric $d$ on $\CalD$ induces a metric on $\CalDq_U$, which we
again call $d$.  Then, $d$ is a $G$-invariant metric on $\CalDq_U$,
and $\CalDq_U$ is complete as a metric space, though neither of these
two facts is obvious, compare \cite[Cor.~3.10]{MS}.

We give $\CalDq_U$ the uniform topology induced by
$U^{}_\mathrm{mACT}$ and again call it the \emph{mixed autocorrelation
  topology}. Again, this is \emph{not} the same as the metric topology
induced by $d$ itself, since it takes small shifts in the sets into
account in order to make the action of $G$ continuous.

\begin{prop} 
   $\CalDq_U$ is complete in the mixed autocorrelation topology.
\end{prop}

\begin{proof}
For $\gL,\gL' \in \CalD_U$, if $\varepsilon >0$ and
$d(\gL,\gL')<\varepsilon$, one finds
\begin{eqnarray*} 
   \varepsilon &>&  d(\gL,\gL') 
   \; = \; \overline{\dens}\,(\gL\,\triangle \, \gL') \\
   &=& \overline{\dens}\left(\big(\gL \backslash (\gL\cap \gL') \big)
      \, \cup \, \big(\gL' \backslash (\gL\cap \gL') \big) \right)\\
   &\ge& \overline{\dens}\,\big(\gL \backslash (\gL \cap \gL' )\big)
   \;\ge \;  \overline{\dens}\,(\gL) - \underline{\dens}\,(\gL \cap \gL') \,.
\end{eqnarray*} 
By symmetry, one also has $\varepsilon > \overline{\dens}\, (\gL') -
\underline{\dens}\, (\gL\cap\gL')$, and hence $| \,
\overline{\dens}\,(\gL) - \overline{\dens}\,(\gL') | < 2\varepsilon$.
 
Now, let $\{\gL_i\} \subset \CalD_U$ be a Cauchy net with respect to
the mixed autocorrelation topology. We want to prove that the net
converges when seen in $\CalDq_U$.
  
{}For any $\varepsilon > 0$ and any open neighbourhood $V$ of $0$,
there is some $n$ so that for all $i,j \succcurlyeq n$ (with
$\succcurlyeq$ referring to the partial order on the index set) and
for suitable $v_{ij} \in V$, one has $d(v_{ij} + \gL_i, \gL_j) <
\varepsilon$. Then, by the above calculation, $|\,
\overline{\dens}\,(v_{ij}+ \gL_i) - \overline{\dens}\,(\gL_j) | <
2\varepsilon$.  Since $\overline{\dens}\,(v_{ij}+ \gL_i) =
\overline{\dens}\,(\gL_i)$, we see that $\{\,
\overline{\dens}\,(\gL_i)\}$ is a Cauchy net in $\RR$ and so converges
to some limit $c\ge 0$. If $c=0$, then $\{\gL_i\} \to \varnothing \in
\CalDq_U$, and we are done. So we only need to consider the case that
$c>0$.

Returning to the Cauchy net $\{\gL_i\} \subset \CalD_U$, choose any
open neighbourhood $V$ of $0$ so that $-V+V+V \subset U$, and any
$\varepsilon$ that satisfies $0<\varepsilon < c/3$.  Fix $n$ so that $i,j
\succcurlyeq n$ implies that $(\gL_i, \gL_j) \in U^{}_\mathrm{mACT}
(V,\varepsilon)$ and $\overline{\dens}\,(\gL_i) > c/2$.

We know that, for all $j,k \succcurlyeq n$,
\[
   d(v_{jk} + \gL_j, \gL_k) < \varepsilon\, , \quad \mbox{for some} \;
   v_{jk} \in V .
\]
Then, for all $j,k \succcurlyeq n$, $d(v_{jk} + v_{nj} + \gL_n, v_{nk}
+ \gL_n) < 3\varepsilon$, or, using translation invariance, $d(-v_{nk}
+ v_{jk} + v_{nj} + \gL_n, \gL_n) < 3\varepsilon$. However, for $x \in
\gL_n$,
\[  
   \{-v_{nk} + v_{jk} + v_{nj}+ x \} \cap \gL_n = 
   \begin{cases} \{x\}, & \text{if} \; -v_{nk} + v_{jk} + v_{nj}=0, \\
   \varnothing , & \text{otherwise} , \end{cases}
\]
because $-V+V +V \subset U$ and all the sets $\gL_{\ell}$ lie in
$\CalD_U$.  So, if $-v_{nk} + v_{jk} + v_{nj}\ne 0$, then $d(-v_{nk} +
v_{jk} + v_{nj} +\gL_n, \gL_n) = 2 \,\overline{\dens}\,(\gL_n) \ge
2(c/2) > 3\varepsilon$, a contradiction. Thus $ v_{jk} + v_{nj}=v_{nk}
$ and $v_{jk} = -v_{kj}$ for all $j, k$.

In principle, $v_{jk}$ depends on $V$ and $\varepsilon$. However, the
same little argument shows that it is actually unique in the sense
that it will be the same element for any $V' \subset V$, $\varepsilon
\le \varepsilon'$.

Let $\gL'_j := -v_{nj} + \gL_j \in \CalD_U$, for all $j \succcurlyeq
n$. Then, for all $j,k \succcurlyeq n$,
\[
   d(\gL'_j, \gL'_k) = d(-v_{nj} + \gL_j, -v_{nk} + \gL_k) 
= d(v_{jk} + \gL_j, \gL_k) < \varepsilon \, . 
\] 
This shows that $\{\gL'_n\}$ is a Cauchy net in $\CalD_U$, with
respect to the metric topology defined by $d$. By \cite[Cor.~3.9]{MS},
it converges to some $\gL \in \CalDq_U$. It is easy to see that also
$\{\gL_n\}$ converges to $\gL$, which completes the argument.
\end{proof}

Denote the equivalence class of $\gL\in\CalD_U$ by $[\gL]$, and
let $\beta$ be the canonical mapping from $\CalD_U$ to 
$\CalDq_U$, i.e., 
\[ 
   \beta : \; \CalD_U \longrightarrow \CalDq_U \, , \quad 
   \varLambda \mapsto [\varLambda] \, .
\]
Each $\gL$ in $\CalD_U$ gives rise to the \emph{autocorrelation hull}
$\bA (\gL)$ defined as the closure of the orbit $G + \beta(\gL)$ in
the mixed autocorrelation topology. By construction, one may as well
consider $\bA (\gL)$ to be the Hausdorff completion of $G$ with
respect to the uniform topology on $G$ that is given by pulling back
the autocorrelation topology from $\CalD_U$. In detail, define a
pseudo-metric (relative to $\gL$) on $G$ by
\begin{equation}\label{MetricOnG}
   d_G (s,t) \; := \; d (t + \gL,s + \gL)
             \;  = \;  d (t -s + \gL, \gL).
\end{equation}
Then, the uniformity on $G$ is described by the sets
\begin{equation}\label{ACTonG1}
   \{(t,s)\in G\times G \, : \, d (v+ t + \gL,s + \gL) < \varepsilon \} 
\end{equation}
where $v\in V$, and $V$ and $\varepsilon$ run over all neighbourhoods
of $0$ and all non-negative real numbers, respectively.

This can be written in a more suggestive way via the set
of $\varepsilon$\emph{-almost periods} of $\gL$,
\begin{equation} \label{eps-periods}
   P_\varepsilon \; := \; 
   \{t \in \gL -\gL : d_G (t,0) <\varepsilon\} \, .
\end{equation}
Then, the entourages (\ref{ACTonG1}) are just the sets 
\begin{equation}\label{ACTonG2}
    \{(t,s)\in G\times G \, : \, t-s \in  V+ P_\varepsilon\} \,.
\end{equation}

These entourages are evidently $G$-invariant.  This has an important
consequence: $\bA (\gL)$, now being the completion of the Abelian
group $G$ with respect to the invariant uniformity as defined by
(\ref{ACTonG2}), carries a natural Abelian group structure.  Moreover,
$G$ acts minimally on $\bA (\gL)$ through the translation action. 
This is the second topological dynamical system for our group $G$.
 Of course, this construction
depends entirely (and crucially) on the starting set $\gL$.  Below, we
shall often shift back and forth between the two views of $\bA (\gL)$:
as a subset of $\CalDq_U$ and as a completion of $G$.

\smallskip 
If we start with a set $\gL \in \CalD_U$, we can form the two hulls
$\beta(\XX (\gL))$ and $\bA (\gL)$. In general, these are \emph{not}
related in any obvious way. In particular, neither is contained in the
other. If, however, $\beta$ is continuous, then obviously $\beta (\XX
(\gL))\subset\bA (\gL)$.  We refer to this mapping $\beta\! : \, \XX
(\gL) \longrightarrow \bA (\gL)$ as the \emph{canonical torus map}.
Moreover, if $\beta$ is continuous and $\XX (\gL)$ is compact, then
$\beta (\XX (\gL)) = \bA (\gL)$, as $\beta (\XX (\gL))$ is then a
compact, and hence closed, set containing $G + \beta (\gL)$.  In this
case, $\bA (\gL)$ becomes a compact topological group.  We shall have
more to say about this situation.

\subsection{Cut and project schemes and model sets}\label{model}

Here, we introduce model sets and discuss some of their basic
features. For further details and proofs, we refer to
\cite{Moody,Martin1,Martin2}.

\smallskip
Model sets arise as (partial) projections from a high dimensional
periodic structure to a lower dimensional subspace. This is formalized
in the following notion.

\smallskip
A \emph{cut and project scheme}, or CPS for short, is a triple
$(G,H,\LL)$ consisting of locally compact Abelian (LCA) groups
$G$ and $H$, with $G$ also being $\sigma$-compact, and a lattice
$\LL$ in $G \times H$ such that the two natural projections
$\pi_1 \! : \, G\times H \longrightarrow G$, $(t,h)\mapsto t$ and
$\pi_2 \! : \, G\times H \longrightarrow H$, $(t,h)\mapsto h$ of the
scheme
\begin{equation} \label{cpScheme}
   \begin{array}{ccccc}
      G & \stackrel{\pi_{1}}{\longleftarrow} & 
        G \times H & \stackrel{\pi_{2}} 
      {\longrightarrow} & H   \\ 
      && \cup \\  && \LL   
   \end{array} 
\end{equation}
satisfy the following properties:
\begin{itemize}
\item The restriction $\pi_1|_{\LL}$ of $\pi_1$ to
      $\LL$ is injective.
\item The image $\pi_2(\LL)$ is dense in $H$.
\end{itemize}
Let $L:=\pi_1 (\LL)$ and $(.)^\star \! : \, L \longrightarrow H$
be the mapping $\pi_2 \circ (\pi_1|_{\LL})^{-1}$. Note that
${}^\star$ is indeed well defined on $L$ and that it can often
be extended to a larger subgroup of $G$ (such as the rational span 
$\QQ L$ in the Euclidean case), but not to all of $G$.

Moreover, as  $\LL$ is  a discrete and co-compact subgroup of
$G\times H$, the quotient
\[
   \TT \; := \; (G\times H) / \LL 
\]
is a compact Abelian group.  In the standard cut and project setting
with Euclidean spaces only, this group is a torus, compare~\cite{BHP}.
There is an obvious action of $G$ on $\TT$ given by
\[
    x + \big((t,h) + \LL\big) \; := \; 
   (x +t, h)+ \LL \, , \quad x\in G \, .
\] 
Then, $(\TT,G)$ is minimal and hence uniquely ergodic as well (as
$\TT$ is a compact Abelian group).

Given a CPS \eqref{cpScheme} and a subset $S
\subset H$, we define $\oplam (S)$ by
\[
    \oplam (S) \; := \; \{x\in L : x^\star \in S\}.
\]
Then, $\oplam (S)$ is relatively dense if the interior of $S$ is
non-empty and it is uniformly discrete if the closure of $S$ is
compact, see \cite{Moody} for details.

A \emph{model set}, associated with the CPS \eqref{cpScheme}, is a
non-empty subset $\gL$ of $G$ of the form
\[
    \gL \; = \; x  + \oplam (y + W),
\]
where $x\in G$, $y\in H$, and $W\subset H$ is a non-empty compact set
with $W = \overline{W^\circ}$. A model set $ \gL = x + \oplam (y + W)$
is called \emph{regular} if $\theta^{}_H (\partial W)=0$.  A (regular)
model set of the above form is called \emph{generic} if $(y +
\partial W) \cap L^\star =\varnothing$.  Any model set is a Delone
set. Namely, it is uniformly discrete as $W$ is compact and relatively
dense as $W$ has nonempty interior.  In fact, they are even Meyer
sets, because $\gL - \gL \subset \oplam (W - W)$ and $W - W$ is
compact, too, and they are thus also FLC sets.  Moreover, a regular
model set has uniform patch frequencies (i.e., the associated
dynamical system is uniquely ergodic) and a generic model set is
repetitive.

\smallskip 
Our prime concern are model sets and their dynamical systems.  It
turns out that the dynamical system associated with the model set
$\oplam (W)$ may contain sets $\gL'$ which are not model sets
themselves with respect to the given CPS.  It is hard to determine
their precise structure in terms of the window. However, under a
condition called irredundancy (see below for more), all of these sets
$\gL'$ satisfy
\begin{equation}\label{IMS} 
   t +  \oplam(c + W^\circ) \; \subset \; \gL' 
   \; \subset \; t+ \oplam(c + W)
\end{equation}
with suitable $t\in G$ and $c \in H$. This suggests to work right from
the start with sets of the form $t+ \oplam(W^\circ) \subset \gL
\subset t+ \oplam(W)$.  This approach is also taken in \cite{LM} in
order to characterize multi-component model sets. We shall call such
sets \emph{inter model sets}, or IMS for short.

The condition we need reads as follows (see \cite{LM} and
Sections~\ref{window} and \ref{toruscp}).

\begin{definition} 
  Let $(G,H,\LL)$ be a CPS. A subset $S$ of $H$ is called
  \emph{irredundant} (with respect to the given CPS), if its
  stabilizer in $H$ is trivial, i.e., if the equation $c + S = S$
  holds only for $c=0\in H$.
\end{definition}

To state our results, we also need the
following definition.

\begin{definition} 
   A dynamical system $(\XX,G)$ is said to be associated with a
   (regular) model set if there exists a (regular) model set
   $\gL$ such that $\XX=\XX(\gL)$.
\end{definition}

\subsection{The torus parametrization:\ Abstract results}
\label{torus} In this paragraph, we look briefly at
  factors of dynamical systems $(\XX,G)$ in which the factors are of
  the form of compact Abelian groups with minimal $G$-actions.  These
  results are essentially known. Throughout, $G$ will be an LCA group
  (although most of this works for other groups as well).  The
  situations that we have in mind are special actions of $G$ by
  translations on point set dynamical systems. These actions
  generalize concepts from \cite{BHP} and \cite{Martin2} known as
  torus parametrizations, and we retain this terminology here.

\begin{definition} 
  Let\/ $\XX$ be a compact space and $(\XX,G)$ a topological dynamical
  system under the action of $G$.  A continuous $G$-map $\rho \!:\,
  \XX \longrightarrow \KK$ into a compact Abelian group $\KK$ on which
  $G$ acts minimally is called\/ \emph{torus parametrization}.
\end{definition}
\begin{definition} \label{nonsingular}
   Let $\rho \!:\, \XX \longrightarrow \KK$ be a torus
   parametrization. For $\xi \in \KK$, we call the inverse image
   $\rho^{-1} (\{\xi\})$ the \emph{fibre} over $\xi$.  Then, $ \gG \in
   \XX $ is called \emph{singular} if the fibre over $\rho (\gG)$
   consists of more than one element.  Otherwise, $\gG$ is called\/
   \emph{non-singular}. In this case, $\{\gG\} = \rho^{-1}(\rho(\gG))$
   is called a \emph{singleton fibre}.
\end{definition}
\begin{lemma} \label{onto}
   If $\rho :\; \XX \longrightarrow \KK$ is a torus
   parametrization,  $\rho$ is onto.
\end{lemma}
\begin{proof} 
As $\XX$ is compact and $\rho$ continuous, the image $\rho(\XX)$ is
compact. Let $\gG$ be an arbitrary element of $\XX$.  As $\rho$ is a
$G$-map, $\rho(\XX)$ contains the orbit of $\rho(\gG)$. As $G$ acts
minimally on $\KK$, this orbit is dense in $\KK$. Thus, $\rho(\XX)$ is
a dense compact subset of $\KK$, hence agrees with $\KK$.
\end{proof}

Let us continue with an interesting property of the torus
parametrization. Namely, each torus parametrization induces a minimal
subsystem of the original dynamical system.

\begin{prop} \label{prop.minimal} 
   Let $\rho \!:\, \XX \longrightarrow \KK$ be a torus
   parametrization. If the set
\[
   \mbox{R}(\XX) \; := \; 
   \{\gG\in \XX : \gG \;\:\mbox{is non-singular}\}
\]
   is non-empty, it is $G$-invariant, and $G$ acts minimally on its
   closure\/ $\XX_R :=\overline{ R (\XX )}$.
\end{prop} 
\begin{proof} The $G$-invariance of $\XX_R$ is 
clear, as $\rho$ is a $G$-map; it remains to show minimality.  To do
so, let an arbitrary $\gG \in R(\XX)$ be given, and consider some
$\gL'\in \XX_R$. Let $\XX (\gL')$ be the closure of its $G$-orbit in
$\XX$. Of course, the restriction $\rho^{}_{\XX (\gL')} \!:\, \XX
(\gL') \longrightarrow \KK$ of $ \rho$ to $\XX ( \gL')$ is a torus
parametrization as well. In particular, it is onto. Thus, we can find
$\gG'\in \XX (\gL')$ with $\rho (\gG') = \rho (\gG)$. By $\gG\in R(\XX
)$, we infer $\gG = \gG'\in \XX (\gL')$. As $\gG\in R(\XX)$ was
arbitrary, this implies $ R(\XX ) \subset \XX (\gL')$. As $ \XX (\gL')
\subset \XX_R$ is clear anyway, we obtain, after taking closures,
\[
     \XX_R \; \subset \; \XX (\gL')
           \; \subset \; \XX_R  \, .
\]
As $\gL'\in \XX_R$ was arbitrary, the statement follows. 
\end{proof}

We now discuss continuity properties of the inverse of a torus
parametrization. While these results are not particularly hard to
prove, they are a crucial ingredient behind the reconstruction of the
window given in Lemma~\ref{cornerstone} in Section~\ref{window}.

\begin{prop} \label{prop.continuity}
   Let $\rho \!:\, \XX \longrightarrow \KK$ be a torus
   parametrization.  Let $\alpha \!:\, \KK\longrightarrow \XX $ be any
   section of $\rho$ $($i.e., $\rho \circ \alpha$ is the identity on
   $\KK\,)$.  Then, $\alpha$ is continuous at all points which are
   images of non-singular points, i.e., at all points of\/
   $\rho(R(\XX ))$.   
\end{prop}


The proof of this proposition is an immediate consequence of the
following lemma.

\begin{lemma} \label{lemma.continuityabstract}  
   Let $K_1$ and $K_2$ be compact spaces and $\sigma \!:\,
   K_1\longrightarrow K_2$ continuous. Let $\xi_1\in K_1$ and
   $\xi_2\in K_2$ be given such that $\{\xi_1\} = \sigma^{-1}
   (\{\xi_2\})$.  Then, a net $(\xi_\iota)$ in $K_1$ converges to
   $\xi_1$ whenever $(\sigma(\xi_\iota))$ converges to $\xi_2$.
\end{lemma}
\begin{proof} 
By compactness of $K_1$, the net $(\xi_\iota)$ has converging
subnets. Thus, it suffices to show that every converging subnet
converges to $\xi_1$. So, consider a converging subnet. Without loss
of generality, we may assume this converging subnet to be
$(\xi_\iota)$ itself. Let $\xi_1'$ be its limit. Then, by continuity
of $\sigma$, we have $\sigma(\xi_1') = \lim_\iota \sigma (\xi_\iota) =
\xi_2$. As, by assumption, $\{\xi_1\} = \sigma^{-1} (\{\xi_2\})$, we
infer $\xi_1' = \xi_1$, and the proof (both of
Lemma~\ref{lemma.continuityabstract} and
Proposition~\ref{prop.continuity}) is complete.
\end{proof}

\section{Outline of the paper and summary of the main theorems}
\label{results}

The overall objective of the paper is to prove Theorem~\ref{main1},
particularly in the direction of sufficiency. The basic setting is
that of a Meyer set $\gL$ for which the local hull $(\XX(\gL),G)$ is
uniquely ergodic. We are interested in continuous $G$-mappings from
the local hull to the autocorrelation hull, and particularly in those
that are non-singular almost everywhere.  Simply the existence of such
a mapping $\ba \! :\, \XX (\gL) \longrightarrow \bA(\gL)$
produces the first prerequisite for the appearance of model sets, a
cut and project scheme.  This is described in Section~\ref{cut}. Any
cut and project scheme (CPS) has associated with it a compact Abelian
group $\TT$ -- the quotient of the product of the ambient group and
the internal group by the associated lattice. A key feature of the cut
and project scheme that arises in our situation is that the mapping
$\ba \! :\, \XX (\gL)\longrightarrow \bA(\gL)$ can be viewed as a mapping
$\XX (\gL)\longrightarrow \TT$:

\smallskip
\textbf{Theorem~\ref{prop.torusparametrization}.}  
\textit{Let $\gL$ be a Meyer set for which $(\XX(\gL),G)$ is uniquely
ergodic. Suppose that there exists a continuous $G$-map $\ba \!:\; \XX
(\gL) \longrightarrow \bA (\gL)$. Then, there is a CPS\/
$(G,H,\LL)$ with associated compact Abelian group $\TT$ for
which $\bA \simeq \TT$ via a topological isomorphism which is a
$G$-map that sends $\gL \in \bA$ to $0 \in \TT$. In particular, there is
a $G$-map $\bt \!:\, \XX (\gL) \longrightarrow \TT$.}

\smallskip 
Having constructed a cut and project set, we next need a window to be
in the regime of model sets.  As studied in Section~\ref{window}, the
crucial condition to provide a window is non-singularity of the
$G$-map $\ba$.  To avoid technical difficulties, we state the result
here in a slightly simplified form.

\smallskip
\textbf{Theorem~\ref{intermediate}A.}  
\textit{Let $\gL$ be a Meyer subset of\/ $G$ such that $(\XX (\gL),G)$ is
uniquely ergodic. Assume that there exists a continuous $G$-map $\ba
\!:\, \XX (\gL)\longrightarrow \bA (\gL)$ which is one-to-one at least at
one point. Then, there is a minimal dynamical subsystem
$(\XX(\gL)_R,G)$ of $(\XX(\gL),G)$ that is associated with a
repetitive model set. In particular, if $\gL$ is repetitive, 
$(\XX(\gL),G)$ itself is associated with a model set. }

\smallskip 
The previous theorem does not assert that the constructed model set is
regular, i.e., that the measure of the boundary of the window is $0$.
Concerning this topic, our result is
Theorem~\ref{prop.measureboundary}.  It shows that the boundary has
Haar measure $0$ if and only if the map $\ba$ is one-to-one almost
everywhere. In fact, if the canonical map $\beta \!:\, \XX
(\gL)\longrightarrow \bA(\gL)$ is one-to-one almost everywhere, we can
get further:

\smallskip
\textbf{Theorem~\ref{main2}.}  
\textit{Let $G$ be a $\sigma$-compact LCA group and $\gL$ a Meyer subset
of $G$ such that the canonical map $\beta \!:\, \XX
(\gL)\longrightarrow \bA (\gL)$ is continuous and one-to-one almost
everywhere, with respect to the Haar measure on\/ $\bA (\gL)
=\beta(\XX (\gL))$. Then, $\XX(\gL)$ is uniquely ergodic and $\gL$
agrees with a regular model set up to a set of density $0$. Furthermore,
if $\gL$ is repetitive, $\XX (\gL)$ is actually associated to a
regular model set.}

\smallskip

So far, we have assumed existence of a continuous $G$-map $\ba\! :\,
\XX(\gL) \longrightarrow \bA(\gL)$.  But what conditions are required
to obtain such a map?  This is studied in Section~\ref{continuity}.
Our main answer is the following.

\smallskip
\textbf{Theorem~\ref{main4}.} 
\textit{Let $\gL$ be a Meyer subset of\/ $G$ such that\/ $(\XX(\gL),G)$
is uniquely ergodic. Then, the following assertions are equivalent.
\begin{itemize}
\item[(a)] There exists a continuous $G$-map $\ba :\; \XX
           (\gL) \longrightarrow \bA(\gL)$.
\item[(b)]  $(\XX (\gL),G)$ has pure point dynamical spectrum with 
           continuous eigenfunctions. 
\end{itemize}
In this case, $\gG,\gG'\in \XX (\gL)$ satisfy $\ba (\gG) =\ba(\gG')$
if and only if $f(\gG)=f(\gG')$ for every eigenfunction $f$.}

\smallskip
The proof of the implication (b) $\Longrightarrow$ (a) of this theorem
requires an intermediate step.  From the assumptions on
$(\XX(\gL),G)$, we create a new dynamical system $(\EE, G)$ by
identifying elements of $\XX$ which are indistinguishable by means of
the continuous eigenfunctions. This new space $\EE$ can be given the
structure of a compact Abelian group. This new group is then shown to
be just $\bA (\gL)$.  This is discussed in Section~\ref{pure} and, in
particular, in Theorem~\ref{EequalA}.

\smallskip  
Theorems~\ref{intermediate}, \ref{prop.measureboundary} and
\ref{main4} establish the sufficiency part of our main
Theorem~\ref{main1} (and most of the necessity too).  This is
discussed in Section~\ref{Proof_of_main1}. The link back is provided
in Section~\ref{toruscp} via the following result.
 
\smallskip
\textbf{Theorem~\ref{main3}.}  
\textit{Let a CPS\/ $(G,H,\LL)$ and a non-empty window $W \subset H$ with $W =
\overline{W^\circ}$ and $\theta_H (\partial W)=0$ be given. If $\gL
\subset G$ satisfies $t+ \oplam(W^\circ) \subset \gL \subset t+
\oplam(W)$ for some $t\in G$, then the canonical map $\beta :\; \XX
(\gL) \longrightarrow \bA (\gL)$ is continuous and one-to-one almost
everywhere.  }

\smallskip
Let us make a short comment here: During the process of proving the
above results, we encounter groups $\bA$ and $\TT$ and maps $\ba$ and
$\bt$ from $\XX(\gL)$ into these groups. We show that these groups are
isomorphic and that, in this sense, $\ba$ and $\bt$ agree. In fact, in
retrospect, we can then even show that these maps agree with the
canonical map $\beta$ introduced above. However, this is not at all
clear at the respective times of appearance and, for this reason, we
carefully distinguish these maps and groups.

\smallskip
Finally, our results also imply an interesting characterization of the 
fully periodic case as discussed  in Section~\ref{crystallographic}:

\begin{definition}
  A set $\gL\subset G$ is called \emph{crystallographic} (or
  fully periodic) if its set of periods 
\[
   \mbox{per}(\gL)  \; := \; \{t\in G : t + \gL = \gL\}
\]
  forms a \emph{lattice}, i.e., a co-compact discrete subgroup of $G$.
\end{definition}

\smallskip
\textbf{Theorem~\ref{main5}.}  \textit{Let $G$ be an LCA group and $\gL$ a
uniformly discrete subset of $G$. Then, the following assertions are 
equivalent.
\begin{itemize}
\item[(i)] $\gL$ is  crystallographic.
\item[(ii)] $\gL$ is Meyer and the map $\beta \!:\; 
            \XX (\gL)\longrightarrow \bA(\gL)$ is 
            continuous and injective. 
\item[(iii)] All of the  following conditions hold:
\begin{itemize}
\item[(1)] All elements of\/ $\XX(\gL)$ are Meyer sets. 
\item[(2)] $(\XX(\gL),G)$ is uniquely ergodic.
\item[(3)] $(\XX(\gL),G)$ has pure point dynamical spectrum 
           with continuous eigenfunctions.
\item[(4)] The eigenfunctions separate all points of\/ $\XX (\gL)$. 
\end{itemize}
\end{itemize}
In this case, $(\XX(\gL),G)$ is also minimal, hence strictly ergodic.  
}

\medskip

The paper revolves around the important concept of Meyer sets.  We
have defined a set $\gL \subset G$ to be Meyer if it is a Delone set
and $\gL - \gL$ is contained in a finite number of translates of
$\gL$. We already noted that this implies that $\gL -\gL$ is also a
Delone set (the important point being that it is uniformly discrete).
For $G= \RR^d$, this is an equivalence, and in fact the most common
definition of a Meyer set is a Delone set whose set of differences is
uniformly discrete. This result is due to Lagarias \cite{Lagarias}. In
the Appendix, we show that the two concepts are equivalent if $G$ is
compactly generated. We also show that, in this case, the requirement
that $\gL -\gL$ be uniformly discrete is equivalent to the apparently
weaker statement that for each compact subset $K$ of $G$, the number
of points of $(t+K) \cap (\gL -\gL)$ is finite and uniformly bounded
as $t$ runs over $G$ (Theorem~\ref{lagariasTheorem}).

\section{Consequences of a continuous $G$-mapping $\beta_{\bA} \! : \;
  \XX (\gL)\longrightarrow \bA (\gL)$: Construction of a
  cut and project scheme} \label{reconstructing}

Let $\gL$ be a Meyer subset of $G$ such that the associated dynamical
system $(\XX (\gL),G)$ is uniquely ergodic. As $\gL$ is Meyer, there
is an open neighbourhood $U$ of $0$ in $G$ so that $\gL$ is
$U$-uniformly discrete, i.e., $\gL\in\mathcal{D}_U$.  Moreover, it
also follows that each element of $\XX(\gL)$ is $U$-uniformly
discrete, too. As discussed in Section~\ref{ac hull}, $\gL$ gives rise
to the autocorrelation hull $\bA$, which is an Abelian group.

In this section, we assume that $\bA (\gL)$ is compact and that there
exists a torus parametrization $\ba \! :\; \XX (\gL) \longrightarrow
\bA (\gL)$. We do \emph{not} assume that the map $\beta_\bA$ is given
by the canonical projection $\beta$.

\smallskip
Our objective in this section is to create a cut and project scheme
out of this torus mapping and to show that $\bA(\gL)$ is
$G$-isomorphic with the torus $\TT$ of the associated cut and project
scheme. Section~\ref{window} then shows how non-singularity of the
torus parametrization can be used to provide and study a window.

\smallskip 
Below, we shall freely use notation from Section \ref{definitions}
and, in particular, Paragraph \ref{model}.

\subsection{Establishing Axioms (A1) -- (A4) of \cite{BM}}
\label{ConditionsofBM}

Let $\gL$ be a Meyer subset of $G$ such that the associated dynamical
system $(\XX (\gL),G)$ is uniquely ergodic. We assume the existence of
a torus parametrization $\ba \! :\; \XX (\gL) \longrightarrow \bA
(\gL)$.
   
In order to create a CPS from this data, we rely on the construction
described in \cite{BM}, based on the Dirac comb $\delta_{\gL}$ of our
point set $\gL$. It is defined by $\delta_\gL := \sum_{x\in \gL}
\delta_x$.  The construction now requires that the four assumptions
(A1), (A2), (A$3^+$) and (A4) of reference \cite{BM} hold for the measure
$\delta_{\gL}$. Let us fix an averaging sequence $\mathcal{A}$ of van
Hove type; the result will not depend on this choice, due to the
unique ergodicity of $(\XX (\gL),G)$.

\smallskip
As $\gL$ is Meyer, the measure $\delta_{\gL} = \sum_{x\in \gL}
\delta_x$ is translation bounded, i.e., for all compact $K\subset G$,
there exists a constant $C_K $ with $\sup_{t\in G} \delta_\gL (t + K)
\leq C_K$. This is just the validity of (A1) for the measure
$\delta_{\gL}$.

As $(\XX(\gL),G)$ is uniquely ergodic, the autocorrelation
\begin{equation}\label{existencegamma}
   \gamma \; := \; 
    \lim_{n\to \infty} \frac{1}{\theta_G (A_n)} 
    \sum_{x,y\in \gG \cap A_n} \delta_{x-y}
\end{equation}
exists for every $\gG \in \XX (\gL)$, does not depend on $\gG$, and
equals $\sum_{x\in \varDelta } \eta(x) \delta_x,$ with $\varDelta =
\gL-\gL$ and a suitable positive definite function $\eta \!:\, G
\longrightarrow \CC$. This is assumption (A2) for $\delta_{\gL}$.

Note that $\eta(0)=\mathrm{dens}(\gL)$, and $\eta (x) =0$ whenever
$x\notin \gL -\gL$. In fact, the function $\eta$ is closely connected
to the metric $d$ described above in (\ref{metric}) and
(\ref{MetricOnG}). More precisely, a direct calculation gives
\begin{equation}\label{metricvseta}
   d(s+ \gL, t + \gL) \; := \; \lim_{n \to\infty} \frac{\mathrm{card}
   \left(\big((s + \gL) \symdiff (t+\gL)\big) \cap A_n \right)}
   {\theta_G(A_n)} \; = \; 2 \big(\eta (0) - \eta (t -s)\big) \, .
\end{equation}
The set $\{x\in G : \eta (x) \neq 0\}$ is clearly a subset of $\varDelta$
and hence uniformly discrete, as $\gL$ is Meyer, and this is
assumption (A$3^+$).

Finally, as $\ba$ is continuous, its image $\bA (\gL)$ is compact. By
\cite{MS}, this implies (see Lemma~\ref{lemma.purepoint} below as
well), that $\gammahat$ is a pure point measure on $\Ghat$. This in
turn means that, for each $\varepsilon>0$, the set of
$\varepsilon$-almost periods defined in (\ref{eps-periods}) is
relatively dense in $G$, compare \cite{BM}. This is assumption (A4).

\smallskip
We close this section by noting that the $\varepsilon$-almost periods
do not depend on $\gL$, but only on $\XX (\gL)$. More precisely, by
uniform existence of the autocorrelation \eqref{existencegamma} and
\eqref{metricvseta}, for every $\gG \in \XX (\gL)$, the identities
\begin{equation}\label{pepsilon}
P_\varepsilon =   \{t \in \gL -\gL : d_G (t,0) <\varepsilon\} \; = \; 
   \{t \in G : d_G (t,0) <\varepsilon\} \; = \; 
   \{t \in \gG -\gG : d (t,0) <\varepsilon\}
\end{equation}
hold  whenever $\varepsilon < 2 \eta (0)$.

\subsection{Creating a cut and project scheme} \label{cut}

Here, we use the method of \cite{BM} to construct a CPS out of
$\gamma$ and $\gL$. This is possible since we have just established the
validity of the necessary conditions (A1), (A2), (A$3^+$) and (A4).

\medskip
Let $L$ be the group generated by the set $\varDelta = \gL - \gL$.
Clearly, the pseudo-metric $d$ discussed in (\ref{MetricOnG}) restricts to
$L$ and gives a pseudo-metric $d_L$ by
\begin{equation}
   d_L (s,t) \; := \; d_G (s,t) \; = \; d(s+ \gL, t + \gL) 
   \; = \; 2 \big(\eta (0) - \eta (t -s) \big) \, ,
\end{equation}
where the last equality follows from \eqref{metricvseta}.  The
topology on $L$ defined by this is again called the
\emph{autocorrelation topology}. It makes $L$ into a topological
group.

A fundamental system of neighbourhoods of $0$ in $L$ is given by the
$P_\varepsilon$, $\varepsilon >0$, defined above in
Eq.~\eqref{eps-periods}.  Let $H$ be the Hausdorff completion of $L$
under the autocorrelation topology and let $\phi \!:\, L
\longrightarrow H$ be the corresponding completion map.  It should be
noted that $\phi$ is not injective in general. In fact, if $\gL$ is a
lattice, one finds $H = \{0\}$.

In any case, let $\LL$ be the subgroup $\{ (t, \phi(t))\mid t \in
L \}$.  Then, this subgroup is a lattice in $G \times H$ and we arrive
at a CPS $(G,H, \LL)$ as shown in
\eqref{cpScheme}.  The pseudo-metric $d_L$ on $L$ induces a
corresponding metric $d_H$ on $H$.  Let $B_\varepsilon^H$ denote the
corresponding open ball of radius $\varepsilon$ in $H$. Then,
\begin{equation}
   P_\varepsilon \; = \;  \phi^{-1}
   \bigl( \phi(L) \cap B^H_\varepsilon \bigr) .
\end{equation}

\begin{prop}
   Let $\gL\subset G$ be a Meyer set such that\/ $(\XX(\gL),G)$ is 
   uniquely ergodic. Then, $\varDelta=\gL-\gL$ is totally bounded\/ 
   $($or precompact\/$)$ in the autocorrelation topology. In particular,
   $\overline{\phi(\varDelta)}$ and $\overline{\phi(\gL)}$ are compact
   subsets of $H$.
\end{prop}
\begin{proof} 
The subsets $P_\varepsilon$, $0 < \varepsilon < 2\eta(0)$, form a
fundamental system of neighbourhoods for $0$ in $L$. Fix one of
them. It is relatively dense in $G$ and hence there is a compact $K$
with $G = P_\varepsilon + K$. Let $s \in \varDelta$ and write $s = t +
k$, with $t\in P_\varepsilon$ and $k\in K$. Then, $s -t \in (\varDelta
- \varDelta) \cap K$ which is a finite set $F$ since $\varDelta$ is a
Meyer set (so, $\varDelta - \varDelta$ is uniformly discrete, see the
Appendix).  Finally, $s = t + s-t \in P_\varepsilon +F$, so $\varDelta
\subset P_\varepsilon +F$, showing that $\varDelta$ is totally
bounded.
\end{proof}

Let $\TT = \TT(\gL) := (G\times H) / \LL$ be the
corresponding compact Abelian quotient group.  There is a natural
action of $G$ on $\TT$, defined by letting $x \in G$ act as $(u,v) +
\LL \mapsto (x +u , v) + \LL \in \TT$ for all $(u,v) \in
G\times H$. This way, $\TT$ becomes a dynamical system for $G$, both
measure theoretically (using the Haar measure $\theta_\TT$) and
topologically. The $G$-orbit of $0 \in \TT$ is dense in $\TT$, as is
every other orbit. The homomorphism $\iota \!: \,G \longrightarrow
\TT$ provided by this orbit is not injective in general: its kernel is
$\ker(\phi) \subset L$, the set of \emph{statistical periods} of
$\gL$. Clearly, $\phi$ plays the role of the $\star$-map, wherefore we
once again write $t^\star$ rather than $\phi(t)$ from now on.

\smallskip
Now, the important fact is that the compact group $\TT$ we have just
constructed agrees with $ \bA(\gL)$ defined in Section~
\ref{ac hull}. More precisely, we have the following result from
\cite{MS}, which follows from the definition of $\bA (\gL)$ and the
characterization of $\TT$ as the completion of $G$ in the so-called
mixed topology given in \cite{BM}. For the convenience of the reader,
we sketch a proof.

\begin{prop}\label{TequalA}  
   Let $\gL\subset G$ be a Meyer set such that\/ $(\XX(\gL),G)$ is
   uniquely ergodic, and let $\ba \! : \; \XX(\gL)\longrightarrow
   \bA (\gL)$ be the corresponding torus parametrization.  Then,
   $\TT\simeq \bA(\gL)$, and this isomorphism is a $G$-map when both
   spaces are given their natural $G$-actions.
\end{prop}
\begin{proof} 
  Let $\alpha \!: \, L \longrightarrow G \times L$ be the diagonal
  map.  Then, $\alpha(L)$ is discrete in $G \times L$ and $(G\times
  L)/\alpha(L)$ becomes a topological group in the usual way.
  Furthermore, $G \simeq (G\times L)/\alpha(L)$ via the canonical
  embedding $ x \mapsto (x,0)+\alpha(L)$, and we provide $G$ with a
  new topology this way, called the \emph{mixed topology}.  There is
  a homomorphism of $(G\times L)/\alpha(L)$ into the compact group $\TT
  = (G\times H)/\LL$ defined by $(x,t)+ \alpha(L) \mapsto (x,t^\star)
  + \LL$.  In \cite{BM}, it is shown that, via this map, $\TT$ is the
  Hausdorff completion of $(G\times L)/\alpha(L)$. Therefore, $\TT$
  may be identified with the Hausdorff completion of $G$ in the mixed
  topology and $\iota(G) \subset \TT$ is the Hausdorff space
  \emph{associated} with $G$. Given this construction of $\TT$, we are
  left with the task to relate the mixed topology to the
  autocorrelation topology.

By definition, a basis for the open neighbourhoods of $0$ in $G$, in the
mixed topology, consists of the sets of the form $V + P_\varepsilon$,
$V$ an open neighbourhood of $0$ in the original topology of $G$,
$\varepsilon >0$ (as these are precisely the sets in $G$ which correspond
to the sets $V \times P_\varepsilon + \alpha(L) \subset (G\times
L)/\alpha(L)$ under our isomorphism).  On the other hand, as discussed
in Section~\ref{definitions}, the autocorrelation completion $\bA$ of
$G$ comes about by supplying $G$ with the uniformity induced from
$\CalD$ which has the sets $U_\mathrm{mACT} (V,\varepsilon) = \{(\gL',
\gL'') : \exists v \in V \mbox{ with } d(v+\gL', \gL'') < \varepsilon
\}$.  The corresponding neighbourhoods of $0$ in $G$ are then 
\[
   U_\mathrm{mACT}(V,\varepsilon) (0) \; =  \;
   \{x \in G : \exists v \in V \mbox{ with }
    d(v + \gL, x+ \gL) < \varepsilon\} \, .
\]
Now, the definition of $P_\varepsilon$ implies
\[ 
      U_\mathrm{mACT} (V,\varepsilon) (0) \; = \; V + P_\varepsilon \, ,
\]
and the proof is complete.
\end{proof}

The key consequence of Proposition~\ref{TequalA} is that our map $\ba
\!:\, \XX (\gL) \longrightarrow \bA (\gL)$ can be interpreted as a
continuous $G$ map $\bt \!:\, \XX (\gL) \longrightarrow \TT$. This
gives
\begin{theorem} \label{prop.torusparametrization}
   Let $\gL$ be a Meyer set for which $(\XX(\gL),G)$ is uniquely
   ergodic.  Suppose that there exists a continuous $G$ map $\ba \!:\,
   \XX (\gL) \longrightarrow \bA (\gL)$. Then, there is a CPS
   $(G,H,\LL)$ with associated compact Abelian group $\TT$ for which
   $\bA \simeq \TT$ via a topological isomorphism which is a $G$-map
   that sends $\gL \in \bA$ to $0 \in \TT$. In particular, there is a
   torus parametrization $ \bt \!:\, \XX (\gL) \longrightarrow \TT$.
   \qed
\end{theorem}

\section{Consequences of the existence of non-singular elements: The window} 
\label{window}

We continue to assume that $\gL$ is Meyer such that the associated
dynamical system $(\XX (\gL),G)$ is uniquely ergodic and that there
exists a torus parametrization, i.e., a continuous $G$-map $\ba \! :\;
\XX (\gL) \longrightarrow \bA (\gL)$. In this section, we investigate
some consequences, first that $\ba$ is non-singular at least at one
element, and second that $\ba$ is non-singular almost everywhere.

\subsection{Existence of a non-singular element} 
   
Assume that $\XX(\gL)$ has at least one non-singular element, see
Definition \ref{nonsingular}.  Thus, we have a dynamical subsystem
$\XX (\gL)_R$ that is the closure of the set of non-singular
elements\/ $R(\XX)$ of\/ $\XX(\gL)$, as defined in
Proposition~$\ref{prop.minimal}$.
  
In the previous section, we have constructed a CPS
from $\gL$ as well as a continuous map $\bt \!:\, \XX (\gL)
\longrightarrow \TT$. In this section, we aim at

\begin{theorem}\label{intermediate}
  Let $\gL\subset G$ be a Meyer set such that $(\XX (\gL),G)$ is
  uniquely ergodic. Assume that there exists a continuous $G$-map $\ba
  \!:\, \XX (\gL)\longrightarrow \bA (\gL)$, which is one-to-one at least
  at one point. Then, there is an irredundant CPS
  $(G,H, \LL)$ associated with $\XX(\gL)$ and a subset $W
  \subset H$, $W = \overline{W^\circ}$ compact, so that every
  non-singular element of $\XX(\gL)$ is of the form
\[
    \gG \; = \;  x + \oplam(-h + W^\circ) 
    \; = \; x + \oplam (-h + W)
\]
  for some $(x,h) \in G\times H$.
  
  Each element of $\XX(\gL)_R$ is repetitive and an IMS
  for the window $W$. If $\gL$ itself is repetitive, 
  one has $\XX(\gL) = \XX(\gL)_R$.
\end{theorem}

The proof requires some preparation. The following lemma is one of the
cornerstones of the present work. It says that the mixed
autocorrelation topology, which is defined by statistical information
at infinity, is actually compatible with the local topology, which is
defined by local information, whenever a certain condition is
met. This condition is that $\gG$ is non-singular relative to $\ba$.
As mentioned above, we always assume in this section that $\gL\subset
G$ is a Meyer set such that $(\XX(\gL),G)$ is uniquely ergodic.
\begin{lemma} \label{cornerstone} 
  Let\/ $\mathcal{A}$ be an averaging sequence for\/ $G$ as introduced
  above, and let  $\gG\in\XX (\gL)$ be non-singular. Given
  any positive integer $M$, there is an\/ $\varepsilon =
  \varepsilon(M) >0$ so that
\[
   t \in P_\varepsilon   \quad \Longrightarrow \quad 
   (t+\gG) \cap A_M  \, = \, \gG \cap A_M \,.
\]
\end{lemma}
\begin{proof}  
  By Proposition~\ref{TequalA}, $\TT\simeq \bA (\gL)$. Now, the
  statement can be concluded from Lemma~\ref{lemma.continuityabstract}
  after noticing that $d (\ba (t + \gG), \ba(\gG)) < \varepsilon$
  whenever $t\in P_\varepsilon$. Namely,
  Lemma~\ref{lemma.continuityabstract} then implies that $t + \gG$ and
  $\gG$ are arbitrarily close in the local topology if $\varepsilon$
  is sufficiently small.  As $\gG$ is Meyer and $P_\varepsilon \subset
  \gG - \gG$ by \eqref{pepsilon}, Fact \ref{snapin} implies that
  $\gG$ and $ t + \gG$ actually agree on arbitrarily large compact
  sets, such as $A_M$, if $\varepsilon > 0$ is chosen accordingly.
\end{proof}

As a consequence of Lemma~\ref{cornerstone}, and extending an argument
used before in \cite{Kola}, we can show that every non-singular
element of $\XX (\gL)$ is a model set:
\begin{prop} \label{prop.window} 
  If $\gG$ is a non-singular element of\/ $\XX (\gL)$ with $0\in \gG$,
  one has $\gG = \oplam(W^\circ)=\oplam (W)$, where $W :=
  \overline{\gG^\star}$ and $W = \overline{W^\circ}$.
\end{prop}
\begin{proof} 
By $0\in \gG$, we have $\gG \subset \gG-\gG\subset \gL-\gL \subset L$.
Now, let $x_0 \in \gG$. Choose a positive integer $M$ so that $x_0 \in
A_M$.  Choose $\varepsilon(M)$ according to Lemma~\ref{cornerstone}.

Let $y \in L$ and suppose that $y_{}^\star \in x_0^\star +
B^H_{\varepsilon(M)}$. Then, $y_{}^\star -x_0^\star \in L^\star \cap
B^H_{\varepsilon(M)} = P_{\varepsilon(M)}^\star$, which implies $y-x_0 \in
P_{\varepsilon(M)}$. Then, $x_0 -y \in P_{\varepsilon(M)}$ and, by
Lemma~\ref{cornerstone}, $(x_0 -y + \gG) \cap A_M = \gG \cap A_M$. This
implies $x_0 -y +u = x_0$ for some $u\in \gG$. Then, $y = u \in \gG$,
so $\gG \supset \oplam(x_0^\star + B^H_{\varepsilon(M)})$ and
$ x_0^\star + B^H_{\varepsilon(M)} \subset W$.
This shows that 
\begin{equation} \label{open}
    \oplam(x_0^\star + B^H_{\varepsilon(x_0)}) \;\subset\; \gG \, , 
    \quad \mbox{for all $x_0 \in \gG$},
\end{equation}
where $\varepsilon(x_0)$ is the $\varepsilon(M)$ of the previous
lemma. Now,
\begin{equation*}
   W := \overline{\gG^\star} = \,
   \overline{\bigcup_{x_0\in\gG} (x_0^\star + B^H_{\varepsilon(x_0)})}
   \;\supset\; \bigcup_{x_0\in\gG} (x_0^\star + B^H_{\varepsilon(x_0)}) 
   =: V.
\end{equation*}
Obviously, $V$ is open and contains $\gG^\star$. Thus,
$\overline{W^\circ} \supset \overline{V} \supset \overline{\gG^\star}
\, = \, W$ and $W = \overline{V}=\overline{W^\circ}$.

\smallskip
By \eqref{open}, $\gG = \oplam (V)$. As $\gG$ belongs to $\XX(\gL)$, a
restriction gives a continuous torus parametrization
\[
    \bt |_\gG :\; \XX (\gG) \subset \XX (\gL) 
    \;\longrightarrow \; \TT .
\]
As $\gG$ is non-singular, the torus parametrization $\bt$ and then
even more the torus parametrization $\bt |_\gG$ is one-to-one at
$\gG$.  

We next show $\partial V \cap L^\star = \varnothing$. If $p\in G$
satisfies $p^\star \in \partial V \cap L^\star$, then, by denseness of
$L^\star$, we can find a net $(t_i) \in L$ with $(t_i^\star )$ in $V
\cap L^\star$ and $t_i^\star\to p^\star$. Without loss of generality,
we may assume that $ p - t_i + \gG = \oplam( p^\star - t_i^\star + V)$
converges to some element $\gG' \in \XX(\gG)$.  Then, $\gG \neq\gG'$
as one contains $p$ and the other does not. On the other hand,
for some $(a,b) \in G\times H$,
\[ 
   \bt |_\gG^{} (\gG)  \; = \;  (a,b) + L^\star    
   \; = \; \lim_i \, (a , b+ p^\star - t_i^\star) + L^\star 
   \; = \; \lim_i \bt |_\gG^{} ( p - t_i  + \gG) 
   \; = \; \bt |_\gG^{} (\gG') \, ,
\] 
contradicting the non-singularity of $\gG$.

By $\partial V \cap L^\star = \varnothing$, we have
\[
     \gG  \; = \; \oplam (V) \; = \; \oplam (V \cup \partial V) 
     \; = \; \oplam(\overline{V}) \; = \; \oplam (W) \, .
\] 
As $\gG = \oplam (V)$ and $V\subset W^\circ$, we infer $\gG = \oplam
(W^\circ)$ as well, and the proof is complete.
\end{proof}

We can use Proposition~\ref{prop.window} to show that the CPS we have
just created is irredundant, and also to determine that each element
in the orbit closure of the non-singular elements, i.e., in
$X(\gL)_R$, is an IMS for some translate of the same window $W$.
 
\begin{prop}\label{beta-prop}
  Let a CPS\/ $(G,H,\LL)$ be given, together
  with a window $W\subset H$ that is non-empty, compact, and satisfies
  $W=\overline{W^\circ}$.  Consider an IMS $\gL$ with
  $\oplam (W^\circ ) \subset \gL \subset \oplam (W)$. With\/
  $\TT=(G\times H)/\LL$ as above, the following assertions are
  equivalent.
\begin{itemize}
\item[(i)] There exists a continuous $G$-map $\bt \!:\,
  \XX(\gL)\longrightarrow \TT$ with $\bt (\gL) = (0,0) + \LL$
\item[(ii)] The window $W$ is irredundant, i.e.,  $W = c+W$ implies
  $c=0$.
\end{itemize}
In this case, $\gL'\in \XX(\gL)$ satisfies $\bt (\gL') = (x,h) +
\LL$ if and only if $x + \oplam(-h + W^\circ) \subset \gL'
\subset x + \oplam (-h + W)$ holds.
\end{prop}

\begin{proof}
The implication (ii) $\Longrightarrow$ (i) follows by the argument given in
\cite{Martin2} to prove the case $\gL = \oplam (W)$ (see \cite{LM} as well).

\smallskip
The implication (i) $\Longrightarrow$ (ii) and the last statement will
be proved together. This will be done in three steps. To this end, let
$\bt \!:\, \XX(\gL)\longrightarrow \TT$ be continuous with $\bt(\gL) =
(0,0) + \LL$, and consider an arbitrary $\gL'\in \XX(\gL)$.

\smallskip \noindent {\sc Step} $1$: 
\textit{$ \bt (\gL') = (x,h) +
\LL$ implies $x + \oplam(-h + W^\circ) \subset \gL' \subset x +
\oplam (-h + W)$.}

\smallskip \noindent Let $(x,h)$ be given with $\bt (\gL') = (x,h) +
\LL$, and let $y\in G$ be chosen so that $0\in \gL'':=-y + \gL'$.  Let
$\{t_n + \gL\}_n$, $t_n \in G$, be a net converging to $\gL''$ in
$\XX(\gL)$.  Without loss of generality, we may assume that $0\in t_n
+ \gL$ for all $n$.  Then, in particular, $t_n\in -\gL$ and therefore
$t_m^\star - t_n^\star \in \, W-W$ for all $n,m$.  As $W -W $ is
compact, we may assume that $\{t_n^\star\}^{}_n \to - k \in H$,
possibly after restricting to a subnet.

Now, $\bt(t_n +\gL) = \iota(t_n) + \bt(\gL)$, where, since $\bt$ is a
$G$-map, $\iota(t_n) = (t_n,0) + \LL = (0, -t_n^\star) +
\LL$, which converges to $(0,k)+\LL$ in $\TT$.  Thus, by
continuity of $\bt$, $\bt(\gL'') = (0,k)+\LL$ and $\bt(\gL') =
\bt (y + \gL'') = (y,k) + \LL$.  As, by assumption, $\bt(\gL') =
(x,h) + \LL$, we infer $(y,k) + \LL= (x,h) +
\LL$. This gives
\begin{equation}\label{hilfs}
     x + \oplam(-h + W^\circ) \; = \; y + \oplam(-k + W^\circ) \quad
    \mbox{and} \quad x + \oplam(-h + W) \; = \; y + \oplam(-k + W).
\end{equation}

Consider an arbitrary $z\in \oplam(-k + W^\circ)$, so that $z^\star +
k \in \, W^\circ$.  Then, for all large $n$, $z_{}^\star - t_n^\star
\in \, W^\circ$ and
\[
   z \in \oplam(t_n^\star +W^\circ) \; = \; t_n +\oplam(W^\circ) \;
   \subset \; t_n + \gL .
\]
Thus, $z\in \gL''$ and $\oplam(-k + W^\circ) \subset \gL''$
follows. Adding $y$, and invoking \eqref{hilfs}, we end up with
\[
    x + \oplam(-h + W^\circ) \; \subset \; \gL'.
\] 
Conversely, if $z \in \gL''$, then $z\in t_n + \gL$ for sufficiently
large $n$, so that $z^\star - t_n^\star \in \, W$ and, in the limit,
$z^\star \in -k + W$, i.e., $z \in \oplam(-k +W)$, which implies $\gL'
\subset y + \oplam(-k + W)$. Again, using \eqref{hilfs}, we obtain
\[
   \gL' \; \subset \; x + \oplam(-h + W).
\]

\smallskip \noindent 
{\sc Step} $2$: 
\textit{$c+ W = W$ implies $c=0$, i.e., condition \textrm{(ii)} holds.}

\smallskip \noindent
Note that $c+ W = W$ implies $c + W^\circ = W^\circ$.  As
$W=\overline{W^\circ}$, the boundary of $W$ is nowhere dense. By the
Baire category theorem, there exists then a $d\in H$ with
\[
   \oplam (d+W^\circ) \; = \; \oplam (d + W).
\]
Moreover, $\bt$ is onto by Lemma~\ref{onto}.
Thus, there exist $\gL', \gL''\in \XX(\gL)$ with
\begin{equation}\label{hilfszwei}
   \bt(\gL')  \; = \; (0,-d) + \LL\, ,\quad 
   \bt(\gL'') \; = \; (0,-d -c) + \LL \, .
\end{equation}
By the result of Step 1, this implies  
\[
   \oplam(d + W^\circ) \;\subset\; \gL' \;\subset\; \oplam (d + W) 
   \quad \mbox{as well as}\quad
   \oplam(d + c+  W^\circ) \;\subset\; \gL'' 
   \;\subset\; \oplam (d + c+  W) \, .
\]
By our choice of $d$, and because we both have $c + W= W$ and $c +
W^\circ = W^\circ$, we can infer $\gL'= \gL''$. This, in turn, implies
$\bt(\gL') = \bt(\gL'')$, and $c=0$ follows from \eqref{hilfszwei}.

\smallskip \noindent 
{\sc Step} $3$: 
\textit{$x + \oplam(-h + W^\circ) \subset
\gL' \subset x + \oplam (-h + W)$ implies $ \bt (\gL') = (x,h) +
\LL$. }

\smallskip \noindent
Let $(y,f)$ with $\bt(\gL') = (y,f) + \LL$ be given. 
By Step 1, we then have
\[
   y + \oplam( - f + W^\circ) \;\subset\; 
     \gL' \;\subset\; y + \oplam(- f + W) \, .
\] 
Adding $-x$ yields
\begin{equation}\label{incleins}
    y - x + \oplam(-f + W^\circ) \;\subset\;
   \gL' - x  \;\subset\;  y - x + \oplam (-f + W) \, .
\end{equation}
On the other hand, the assumption on $(x,h)$ gives
\begin{equation}\label{inclzwei}
     \oplam(-h + W^\circ) \;\subset\; \gL' - x 
     \;\subset\; \oplam (-h + W) \, .
\end{equation}
These inclusions show that $(y-x)$ belongs to $L$ and we can rewrite
\eqref{incleins} as
\begin{equation}\label{incldrei}
     \oplam((y-x)^\star  -f + W^\circ) \;\subset\; 
     \gL' - x \;\subset\;  \oplam ( (y-x)^\star - f + W) \, .
\end{equation}
Now, a combination of  \eqref{inclzwei} and \eqref{incldrei} gives
\[ 
   \oplam((y-x)^\star  -f + W^\circ) \;\subset\; 
   \oplam (-h + W) \quad\mbox{and}\quad
   \oplam(-h + W^\circ) \;\subset\; \oplam ( (y-x)^\star - f + W)\, ,
\]
which in turn implies
\[ 
   ((y-x)^\star - f + W^\circ) \cap L^\star 
   \;\subset\; -h + W
   \quad\mbox{and}\quad (-h + W^\circ) \cap L^\star 
   \;\subset\; (y-x)^\star - f + W  \, .
\]
Taking closures and using $\overline{W^\circ} = W$ as well as the
denseness of $L^\star$ in $H$, we obtain
\[ 
   (y-x)^\star - f +  W \subset -h + W 
   \quad\mbox{and}\quad -h + W  \subset (y-x)^\star - f + W.
\]
These inclusions yield  $f - h - (y-x)^\star + W = W$ and, by Step 2, 
\[ 
   f - h - (y-x)^\star \; = \; 0.
\]
This, however, means $(y,f) + \LL = (x,h) + \LL=\bt(\gL')$,
and the proof of Step 3, and also of the entire claim, is complete.
\end{proof}

\subsection{The proof of Theorem~\ref{intermediate} } 
\label{proofintermediate} 

To prove Theorem ~\ref{intermediate}, consider a non-singular element
$\gG$ of $\XX(\gL)$. By translating $\gG$, we may assume $0\in \gG$
without loss of generality.  Proposition~\ref{prop.window} then
implies $\gG = \oplam(W^\circ)= \oplam (W)$, where $W :=
\overline{\gG^\star}$ and $W = \overline{W^\circ}$ is compact. By
Proposition~\ref{TequalA}, $\bA (\gL) \simeq \TT$.

Assume $\gG = \gL$ for the moment. Then, by Proposition
\ref{beta-prop}, every $\gL' \in \XX(\gL)$ is an IMS of
the form that we require. If, on the other hand, $\gL$ is singular,
these results apply to all the elements of $\XX(\gL)_R$, since it
contains all the non-singular elements and is the closed hull of any
of its elements. As pointed out in Fact 2, the elements of
$\XX(\gL)_R$ are all repetitive.

This finishes the proof of Theorem~\ref{intermediate}. \qed

\begin{remark}
There is very little that one can say about the generator $\gL$ of the
hull $\XX(\gL)$ being a model set, or even an IMS, if repetitivity or
some other consistency property is not assumed. One can, for instance,
take a model set, add some finite set of spurious points, and take the
hull of the resulting set. That destroys the set as a model set, but
does not destroy the properties of the minimal part $\XX(\gL)_R$ of
the hull, which will not have been altered.  However, with the
assumption of non-singularity almost everywhere, we can obtain
information up to sets of density $0$.
\end{remark}

\subsection{Consequences of non-singularity almost everywhere}

In this section, we shall prove Theorem~\ref{main2}. To do so, we need
some preparation around the regularity of the window in the cut and
project scheme. To do so, we assume the following setting.
\begin{itemize} \label{S}
\item[\textbf{(S)}] $(G,H,\LL)$ is a CPS,
  $W\subset H$ is a non-empty, compact set with $W=\overline{W^\circ}$,
  and $\gL$ is an arbitrary IMS for it, i.e., $\oplam(W^\circ)\subset
  \gL\subset\oplam( W)$.  There exists a continuous $G$-map $\bt \!:\,
  \XX(\gL)\longrightarrow \TT$ with $\bt(\gL) = (0,0) + \LL$.
\end{itemize}

\smallskip \noindent  
Proposition~\ref{beta-prop} has the following consequence.
\begin{prop}  \label{help}
   Let\/ {\rm \textbf{(S)}} be valid, with an IMS $\gL$. For any
  $c\in H$, the following properties are equivalent.
\begin{itemize}
   \item[(i)] $\oplam (-c+W^\circ) = \oplam(-c+W)$;
   \item[(ii)] $\partial (-c+W) \cap L^\star = \varnothing$;
   \item[(iii)] The fibre over $(0,c) +\LL$ is non-singular.
\end{itemize}
In this case, $\oplam (-c+W^\circ) = \oplam(-c+W)$ constitutes the
fibre over $(0,c) +\LL$, and one has the inclusion\/ $\XX(\oplam
(-c+W)) \subset \XX(\gL)$.
\end{prop}
\begin{proof} 
The equivalence of (i) and (ii) is obvious. Also, \textbf{(S)} allows us
to use Proposition~\ref{beta-prop}, whence we see that (i) implies (iii).

It remains to show that (iii) implies (ii), or its contraposition. 
To this end, let us assume that there is some $p\in
L$ with $p^\star \in \partial(-c + W)$. The $\bt$-fibre over $(0,c)
+ \LL$ is non-empty and consists of the elements $\gL' \in
\XX(\gL)$ such that
\[ 
     \oplam (-c+W^\circ) \subset \gL' \subset  \oplam (-c+W) 
\]
by Proposition~\ref{beta-prop}. We claim that there are at least two
elements on this fibre, one of which contains $p$ while the other does
not.

Take any $\gL'$ on the fibre. Suppose first that $p \notin \gL'$.
Since $p^\star$ is on the boundary of $-c+W$, there is a net
$\{\ell_n \}$ in $L $ with $\{\ell_n^\star \} \longrightarrow c$, such
that $p \in \oplam(-\ell_n^\star + W^\circ)$ for all $n$. Then, on the
fibre over $(0, \ell_n^\star) + \LL$, there is a set $\gL_n \in
\XX(\gL)$ with $\oplam (-\ell_n^\star+W^\circ) \subset \gL_n \subset \oplam
(-\ell_n^\star+W)$. By the compactness of $\XX(\gL)$, there is a convergent
subnet of $\{\gL_n\}$ which we may assume to be $\{\gL_n\}$
itself. Let $\{\gL_n\} \longrightarrow \gL'' \in \XX(\gL)$. Then, $p
\in \gL_n$ for all $n$ implies $p \in \gL''$. Also, $\bt(\gL'') =
\lim_n \, \bt(\gL_n) = \lim_n \, (0,\ell_n^\star) + \LL = (0,c)
+\LL$, so $\gL''$ is on the same fibre as $\gL'$, but it contains
$p$ whereas $\gL'$ does not.

The argument for the case when $p \in \gL'$ is similar. This time,
choose a net $\{\ell_n\}$ in $L$ with $\{\ell_n^\star \} \to c$, $p
\notin \oplam(-\ell_n^\star + W)$, for all $n$.  We then find $\gL_n$
on the fibre over $(0, \ell_n^\star) + \LL$, with $p \notin
\gL_n$, and get $\gL'' \in \XX(\gL)$ on the fibre over $(0,c) + \tilde
L$, also with $p \notin \gL''$.

\smallskip
The last statement of the Proposition is obvious.
\end{proof}

Next, let us relate the properties of $W$ versus $\partial W$ to the
injectivity of $\bt$.

\begin{theorem}  \cite{Moody2001} \label{uniformdistribution}
   Let $(G,H,\LL)$ be a CPS. Let $M$ be a
   measurable, relatively compact set in $H$. Then,
\[ 
    \dens (x + \oplam ( M -h ) ) \; := \;
    \lim_{n\to \infty} \frac{ \mathrm{card}\big( (x + \oplam ( M -h )) 
    \cap A_n \big)} {\theta_G (A_n)} \; = \; 
    \dens (\LL)\,\theta_H (M) ,
\]
   which is valid
   for all\/ $(x,h) \in G\times H$ if\/ $ \theta_H (\partial M)=0$,
   and otherwise for\/ $\theta_G\times \theta_H$-almost every\/ $(x,h)\in G 
   \times H$.  \qed
\end{theorem}

\begin{lemma} \label{WeilResult} 
   Let $M \subset \TT$ be any measurable subset whose preimage in
   $G\times H$ is contained in a subset of the form $G \times B$ with
   $\theta_H (B) = 0$. Then, $\theta_{\TT}(M) = 0$.
\end{lemma}
\begin{proof} 
Observe first that $\theta_{G\times H} (A\times B) = \theta_G (A) \,
\theta_H (B) = 0$ for any relatively compact measurable set $A\subset
G$.  Since $G$ is $\sigma$-compact, we may now employ the averaging
sequence $\mathcal{A}=\{A_n\}$ of Section~\ref{basic}, with
$A_{n}\subset A_{n+1}$ and $G=\bigcup_n A_n$, to conclude that also
$\theta_{G\times H} (G\times B) = 0$.

Let now  $\pi \!: \, G\times H \longrightarrow \TT$ be the canonical
projection.  Define, for $\xi \in \TT$, the measure $\nu_\xi$ on
$G\times H$ by
\[ 
   \nu_\xi  \; := \;
   \sum_{z\in \pi^{-1}(\xi)} \delta_z \, .
\]
Standard desintegration (e.g., using a fundamental domain) shows
that $\theta_{G\times H} = \theta_\TT\circ \nu$, i.e., $\int f (z)
\dd \theta_{G \times H } (z) = \int \nu_\xi (f) \dd \theta_\TT (\xi)$ for
any measurable nonnegative $f$ on $G\times H$. This gives
\[
   0\leq \theta^{}_\TT (M) = \int 1^{}_M (\xi) \dd \theta^{}_\TT (\xi) 
    \leq \int \nu^{}_\xi (1^{}_M \circ \pi ) \dd \theta^{}_\TT (\xi) 
    \leq \int \nu^{}_\xi (1^{}_{G\times B}) \dd \theta^{}_\TT (\xi) = 
    \theta^{}_{G\times H} (G\times B)
\]
and the proof is finished because the last term vanishes as shown above.
\end{proof}

\begin{theorem}\label{prop.measureboundary} 
   Let\/ {\rm \textbf{(S)}} be in place. Then, the boundary of\/ $W$ 
   has measure\/ $0$ if and only if $\bt$ is one-to-one almost everywhere.
\end{theorem}
\begin{proof}  By Proposition~\ref{help}, 
$\gL' \in \XX(\gL)$ is non-singular if and only if $\gL' = x +
\oplam(-h+ W^\circ) = x + \oplam(-h+ W)$ and $L^\star \cap (-h+
\partial W) = \varnothing$ for some $(x,h) \in G\times H$. In this
case, one has $x + \oplam(-h + \partial W) = \varnothing$.  We thus
have $\dens (x + \oplam ( -h +\partial W ) ) = 0$ at this point.

If $\bt \!:\, \XX(\gL) \longrightarrow \TT$ is one-to-one $\TT$-a.e., we
also have this relation $G\times H$-a.e., due to $\theta_{G\times H} =
\theta_\TT\circ \nu$ (see the proof of the previous lemma).
Consequently, by Theorem~\ref{uniformdistribution} and because
$\dens (\LL)\neq 0$, we may conclude that
$\theta_H (\partial W)= 0$.

Conversely, suppose that $ \theta_H (\partial W)= 0$.  By
  Proposition~\ref{help},
\[ 
  \begin{split}
  F & := \, \{\xi\in \TT : \mbox{the fibre over $\xi$ contains 
          more than one element}\} \\
    & \, = \, \{ (x,h)+\LL \in \TT : 
   \oplam(-h + W^\circ) \neq \oplam(-h + W)\} \,.
   \end{split}
\]

 This gives
\begin{eqnarray*} 
  \theta_\TT (F) &= &\theta_{\TT} \big(\{ (x,h)+\LL \in \TT : 
  \oplam(-h + W^\circ) \neq \oplam(-h + W)\}\big) \\
  &=& \theta_{\TT} \big(\{ (x,h) +\LL \in \TT : h \in
 L^\star + \partial W \} \big) \\
  &=& \theta_{\TT}\big(G \times (L^\star + \partial W)\mbox{ mod } 
  \LL \big) \; = \; 0,  \, 
\end{eqnarray*}
where we used Lemma~\ref{WeilResult} in the last step together with
the fact that $L^\star$ is countable.
\end{proof}

We can now proceed to the final result of this section.

\begin{theorem}\label{main2}
  Let $G$ be a $\sigma$-compact LCA group and $\gL$ a Meyer subset of
  $G$ such that the canonical map $\beta \!:\, \XX
  (\gL)\longrightarrow \bA (\gL)$ is continuous and one-to-one almost
  everywhere, with respect to the Haar measure on\/ $\bA (\gL)
  =\beta(\XX (\gL))$.  Then, $\XX(\gL)$ is uniquely ergodic and $\gL$
  agrees with a regular model set up to a set of density $0$.
  Furthermore, if $\gL$ is repetitive, $\XX (\gL)$ is actually
  associated to a regular model set.
\end{theorem}

\begin{proof} 
We are given a Meyer set $\gL$ and assume that $\beta \!:\, \XX
(\gL)\longrightarrow \bA (\gL)$ is continuous and one-to-one almost
everywhere. We first want to show that $\gL$ differs from a model set
up to a set of points of density $0$.  As $\beta$ is continuous, $\bA
(\gL)$ is compact. Moreover, $G$ acts minimally on $\bA (\gL)$ by
definition.  Thus, $(\bA ( \gL),G)$ is uniquely ergodic with the Haar
measure $\theta_\bA$ on $ \bA ( \gL)$. We show that $(\XX (\gL),G)$ is
uniquely ergodic as well.

As $\beta$ is one-to-one almost everywhere, there exists a
subset $\bA'\subset \bA(\gL)$ of full measure such that $\beta$ is
one-to-one on $\XX':=\beta^{-1}(\bA')$ and the complement of $\XX'$ is
mapped into the complement of $\bA'$ by $\beta$. Now, note that, by
Proposition~\ref{prop.continuity}, any inverse of $\beta$ is
continuous when restricted to $\bA'$. Thus, extending this continuous
function, say by setting it constant on $\bA(\gL)\setminus \bA'$, we
find a measurable $\alpha \!:\, \bA\longrightarrow \XX (\gL)$, which
is an inverse to $\beta$ on $ \bA'$.

Let $\mu$ be any $G$-invariant probability measure on $\XX (\gL)$. As
$(\bA (\gL), G)$ is uniquely ergodic, $\beta^\ast (\mu)$ is the Haar
measure $\theta_\bA$ on $\bA(\gL)$. In particular, $\mu
(\beta^{-1}(M)) =0 $ whenever $M$ is a subset of $\bA(\gL)$ of measure
$0$.  In particular, $\mu (\XX(\gL)\setminus \XX') =0$.  Let $f$ be
any measurable bounded function on $\XX(\gL)$. Then, $f$ and $f \circ
\alpha \circ \beta$ only differ on $ \beta^{-1}( \bA(\gL)\setminus
\bA')$, which has $\mu$-measure $0$.  This implies
\[
    \mu (f) \; = \; \mu (f \circ \alpha \circ \beta ) 
    \; = \; \beta^\ast(\mu) (f \circ \alpha) \; = \; 
    \theta_\bA ( f \circ \alpha) \; = \; 
    \alpha^\ast (\theta_\bA) (f) \, .
\]
Thus, $\mu$ is uniquely determined and the unique ergodicity of 
$(\XX(\gL),G)$ follows. 

\smallskip
Now, the assumptions of Theorems~\ref{prop.torusparametrization} and
\ref{intermediate} are satisfied, and we find both a CPS such that
$\bA (\gL) \simeq \TT$ and a dynamical system $(\XX(\gL)_R,G)$
associated to a model set $\gG =\oplam (W)$ inside of $(\XX ( \gL),G)$
with irredundant $W$ and (metrizable) internal group $H$.

{}From the previous results, and Theorem~\ref{prop.measureboundary} in
particular, we know that the almost one-to-one-ness of $\beta$ forces the
boundary of $W$ to have measure $0$. Consider the fibre lying over
$(x,h) + \LL$. If the fibre is non-singular, the single element of
$\XX(\gL)$ is $x + \oplam(h+W)$, which is a \emph{regular} model set.
Even if the fibre is singular, set(s) lying there differ by density
$0$ from the regular model set $x + \oplam(h+W)$, since $\dens( x +
\oplam(h +\partial W)) =0$ by Theorem \ref{uniformdistribution}.

Of course, if $\gL$ is repetitive, $\XX(\gL)$ is generated by any of
its elements, and so $\XX(\gL)$ is actually associated to a regular
model set.
\end{proof}

\section{Existence of a continuous $\ba$ and pure point spectrum with
  continuous eigenfunctions }
\label{continuity}

Let $\gL$ be a Meyer set such that $(\XX(\gL),G)$ is a uniquely
ergodic dynamical system. The existence of a torus parametrization
$\ba \!:\, \XX(\gL) \longrightarrow \bA(\gL)$ has proved to be the key
to linking $\gL$ to the realm of model sets. In this section, we connect
the existence of a torus parametrization with properties of the
dynamical system $\XX(\gL)$ itself. These properties are pure
pointedness of the spectrum and continuity of the eigenfunctions. 

\begin{theorem}\label{main4} 
   Let $\gL$ be a Meyer subset of\/ $G$ such that\/
   $(\XX(\gL),G)$ is uniquely ergodic. Then, the following
   assertions are equivalent.
\begin{itemize}
\item[(a)] There exists a continuous $G$-map $\ba :\; \XX
           (\gL) \longrightarrow \bA(\gL)$.
\item[(b)]  $(\XX (\gL),G)$ has pure point dynamical spectrum with 
           continuous eigenfunctions. 
\end{itemize}
In this case, $\gG,\gG'\in \XX (\gL)$ satisfy
$\ba (\gG) =\ba(\gG')$ if and only
if $f(\gG)=f(\gG')$ for every eigenfunction $f$.
\end{theorem}

This and the next section of the paper are devoted to the proof of
this result.  In this section, we prove Theorem~\ref{main4} in the
direction (a) $\Rightarrow$ (b). In the following section, we prove the
converse.

\subsection{The proof of (a) $\Longrightarrow$ (b) of Theorem~\ref{main4} }
Let  $\gL$ be a Meyer set  with associated
uniquely ergodic dynamical system $(\XX (\gL),G,\mu)$ and let $T$
be the corresponding unitary representation of $G$ on
$L^2 (\XX (\gL),\mu)$. 

Recall that the eigenvalues of this dynamical system form a subgroup
of $\widehat{G}$, which we denote by $P(T)$. In addition, we need to
consider the diffraction measure $\widehat{\gamma}$, which is constant on
$\XX(\gL)$ due to the unique ergodicity. For any measure $\nu$ on
$\widehat{G}$, we introduce the set
\begin{equation} \label{defOfP}
    P(\nu) \, := \, \{ k \in \widehat G \, : \, \nu(\{k\}) \ne 0 \},
\end{equation}
which is a countable set, and the subgroup of $\widehat G$ that it
generates, denoted by $\langle P(\nu) \rangle$.

We recall the following result that has already been established in
the literature.
\begin{lemma} \label{lemma.purepoint} 
   Let $\gL\subset G$ be a Meyer set.
   If\/ $(\XX (\gL),G)$ is uniquely ergodic, the following
   assertions are equivalent.
\begin{itemize}
 \item[(i)] $\bA(\gL)$ is compact.
 \item[(ii)] $\gammahat$ is a pure point measure.
 \item[(iii)] $(\XX(\gL),G)$ has pure point dynamical spectrum. 
\end{itemize}
   In this case, the dynamical spectrum\/ $P(T)$ of\/
   $(\XX(\gL),G)$ satisfies\/ $P(T) = \langle P(\gammahat) \rangle $. 
\end{lemma}
\begin{proof}
The equivalence of (i) and (ii) is shown in \cite{MS}.  The
equivalence of (ii) and (iii) is proved in \cite[Thm.\ 3.2]{LMS-1}. 
The last statement is proved in \cite[Thm.\ 9]{BL}.
\end{proof}

If there exists a continuous $G$-map $\ba :\; \XX
           (\gL) \longrightarrow \bA(\gL)$, 
Theorem~\ref{prop.torusparametrization} tells us that we have a CPS
$(G,H,\LL)$ and a compact group $\TT=(G\times H)/\LL$. Moreover,
$\bA (\gL)=\TT$.  Thus $\ba$ induces a continuous map $\bt$ between
$\XX (\gL)$ and $\TT$.  There is then a canonical homomorphism $ \iota
\!:\, G \longrightarrow \TT$ of topological groups with dense range
defined by $x \mapsto (x,0) + \LL$. Dualizing, we obtain an injective
homomorphism $\hat{\iota} \!:\, \widehat{\TT}\longrightarrow \Ghat$ of
the dual topological groups. Lemma \ref{lemma.purepoint} tells us that
$\gammahat$ is a pure point measure.

\begin{lemma} \label{lemma.support}
    Let $\gL$ be a Meyer set such that\/ $(\XX (\gL),G)$
    is uniquely ergodic and \/ $\gammahat$ is a pure point
    measure. Then, $\langle P(\gammahat) \rangle \subset \hat{\iota}
    (\widehat{\TT})$.
\end{lemma}

\begin{proof}
  Due to unique ergodicity, each element of $\XX(\gL)$ has the same
  autocorrelation measure $\gamma$. Let $C_\mathrm{c}(G)$ denote the
  space of continuous complex-valued functions of compact support on
  $G$.  For every $c \in C_\mathrm{c}(G)$, we define $g_c \!:\, G
  \longrightarrow \CC$ by $g_c = c\ast \tilde{c} \ast \gamma$. 
  Then, there is a continuous positive definite function
  $g_c^{\TT} $ on $\TT$ so that $g_c^{\TT}\circ \iota = g_c$ 
  (see Section~4 of \cite{BM} as well).  In
  particular, we can expand $g_c^{\TT} $ in a uniformly converging
  Fourier series
\[
      g_c^{\TT}  (x) \; = \; \sum_{k  \in \widehat{ \TT}} a_c (k)
       (k, x)
\]
with nonnegative numbers $a_c (k)$ that satisfy
\[
      \sum_{k  \in \widehat{ \TT}}  a_c (k)
       \; = \; g_c (0) \, .
\]
Composing $ g_c^{\TT} $ with the homomorphism $\iota \!:\, G
\longrightarrow \TT$ and using the definition of $\hat{\iota}$, we
obtain
\[
     g_c (x) \; = \; \sum_{k \in \widehat{\TT}} a_c (k)\,
      (\hat{\iota}(k), x) \, .
\]
As the $a_c (k)$ are summable, we can calculate the
Fourier transform of $g_c$ to arrive at
\[
     |\widehat{c}|^2 \gammahat  \; = \;
      \widehat{g_c} \; = \; \sum_{k \in \widehat{\TT}}
       a_c (k) \, \delta_{\hat{\iota}(k)} \, ,
\]
which is a finite positive measure on $\widehat{G}$.

This shows that
\[
     B:= \{ k\in \widehat{\TT} : a_c (k) >0 \;\mbox{for
          some continuous $c$ with compact support} \}
\]
is mapped unter $\hat{\iota}$ into $P(\gammahat)$ defined in
\eqref{defOfP}.  Taking for $c$ an approximate unit, one infers that
$B$ is actually mapped onto $P(\gammahat)$.  As $\hat{\iota}
(\widehat{\TT})$ is a subgroup of $\Ghat$, which then contains
$\iota(B)$, the desired conclusion follows immediately.
\end{proof}

\begin{prop} \label{betaimpliespp}
    Let $\gL$ be a Meyer set in $G$ such that $(\XX (\gL),G)$ is
    uniquely ergodic. If there exists a continuous $G$-mapping $\ba \!
    : \; \XX (\gL)\longrightarrow \bA (\gL)$, then $(\XX (\gL),G)$ has
    pure point dynamical spectrum with continuous eigenfunctions.
\end{prop}

\begin{proof} 
By Lemma~\ref{lemma.purepoint}, the dynamical system $(\XX (\gL),G)$
has pure point dynamical spectrum. Moreover, as discussed after the
lemma, we can then identify $\bA (\gL)$ and $\TT$. Thus, $\ba $ yields 
a continuous $G$-map from $\XX(\gL)$ to $\TT$.

Every element $\lambda\in \widehat{\TT}$ gives rise to a
continuous eigenfunction
\[
     f_\lambda := \lambda \circ \bt : \;
     \XX (\gL)\longrightarrow \CC
\]
to the eigenvalue $\hat{\iota} (\lambda)$, and we infer
\[
      \hat{\iota}(\widehat{\TT})\subset P(T) \, ,
\]
where the point spectrum $P(T)$ is the set of eigenvalues. Combining
this with the results of Lemma~\ref{lemma.purepoint} and
Lemma~\ref{lemma.support}, we obtain the following chain of
inclusions:
\[
     \hat{\iota}(\widehat{\TT})\; \subset \;
      P(T) = \langle P(\gammahat)\rangle  \; \subset \;
     \hat{\iota} (\widehat{\TT}) \, .
\]
Therefore, $\hat{\iota}(\widehat{\TT}) = P(T)$.  Thus, the
$f_\lambda$, $\lambda\in \widehat{\TT}$, provide eigenfunctions for
each eigenvalue. As each eigenvalue has multiplicity one by
ergodicity, we have found a complete system of eigenfunctions, all of
which are continuous.
\end{proof}

\noindent
This finishes the proof of Theorem~\ref{main4} in the direction
$(a)\Rightarrow (b)$ . \qed

\begin{remark}
Under the hypotheses of Proposition~\ref{betaimpliespp},
$\langle P(\gammahat) \rangle = \hat{\iota} (\widehat{\TT})$.
This also holds under the assumptions of Lemma~\ref{lemma.support} 
(and, in fact, in much more general situations as well).
This will be discussed further in \cite{BLM}.
\end{remark}

\section{Consequences of unique ergodicity, pure point dynamical
spectrum, and continuous eigenfunctions}\label{pure}

The aim of this section is to prove the following theorem. As
discussed at the end of this section, this theorem will provide the
proof of the missing direction of Theorem~\ref{main4}.

\begin{theorem}  \label{EequalA}  
  Let $\gL\subset G$ be a Meyer set and\/ $(\XX (\gL),G)$ be uniquely
  ergodic. Suppose that\/ $(\XX(\gL),G)$ has pure point dynamical
  spectrum and continuous eigenfunctions. Let $\bA (\gL)$ be the
  autocorrelation hull. Then, there exists a torus parametrization
  from $\XX$ to $\bA$ for which $\gL \mapsto 0$.
\end{theorem}

To prove this result, we proceed as follows.  In
Paragraph~\ref{general}, we assume that $(\XX,G)$ is an arbitrary
uniquely ergodic dynamical system with pure point spectrum and
continuous eigenfunctions. We then show how to construct a compact
topological group $\EE$ and a continuous surjective $G$-map $\be \!:\,
\XX \longrightarrow \EE$.  In the subsequent paragraphs, we return to
the case of $(\XX,G)$ being a Meyer dynamical system, assuming now
that we have pure point spectrum and continuous eigenfunctions, and
show that the continuous map $\be \!:\, \XX \longrightarrow \EE$
constructed in the first paragraph is effectively none other than a
torus parametrization $\ba \!:\, \XX \longrightarrow \bA$.

\begin{remark} 
 As investigated by Robinson \cite{Rob} in the case of $G=\RR^d$ and
  $G=\ZZ^d$, continuity of the eigenfunctions is related to uniform
  existence of certain limits (see \cite{Len} for recent results in
  the case of general LCA groups $G$ as well).
\end{remark}

\subsection{A general construction}\label{general}

Let $(\XX,G)$ be a uniquely ergodic dynamical system with unique
$G$-invariant probability measure $\mu$.  This gives rise to a unitary
representation $T$ of $G$ on $L^2 (\XX,\mu)$.  

Assume that $T$ has pure point dynamical spectrum with all
eigenfunctions continuous. This means that $L^2(\XX,\mu)$ has an
orthonormal basis $\{f_\lambda$ : $\lambda \in P(T)\}$ where the point
spectrum $P(T)$ (i.e., the set of eigenvalues of $T$) is some
\emph{subgroup} of $\widehat{G}$, the character group of $G$, and each
$f_\lambda$ is a continuous eigenfunction for the character $\lambda$.
Note that, due to ergodicity, all eigenvalues are simple, and the
corresponding eigenspaces are thus one-dimensional \cite{Wal}.  We may
assume that each $f_\lambda$ is normalized to $1$ (in the $L^2$-norm).

Define $\gL' \sim \gL''$ when $f_\lambda(\gL') =
f_\lambda(\gL'')$ for all $\lambda \in P(T)$. Let $\EE := \XX
/\!\!\sim\,$\/ and let $\be$ denote the canonical mapping from $\XX$
to $\EE$. Note that the $f_\lambda$ can be factored through the
equivalence relation. Give the quotient space the uniform structure
for which the cylinder sets given by
\[
  U(F,\varepsilon) \; := \; 
  \{(\be(\gL'),\be(\gL''))\in \EE \times
  \EE : \,|f_\lambda(\gL') - f_\lambda(\gL'')| <\varepsilon,
  \lambda \in F\}\, ,
\]
where $F$ runs through all finite subsets of $P(T)$ and $\varepsilon$
through the positive reals, are a fundamental system of
entourages. The mapping $\be \!:\, \XX \longrightarrow \EE$ is
uniformly continuous because the eigenfunctions are continuous (hence
uniformly continuous, since $\XX$ is compact). Thus, $\EE$ is compact
and hence complete.

Each of the basic entourages of $\EE$ is actually $G$-invariant (since
the $f_\lambda$ are \emph{eigenfunctions}) and we obtain from this an
obvious $G$-action on $\EE$ for which the natural mapping
$\be$ from $\XX$ to $\EE$ is a $G$-map. This implies the
orbit $\be(G + \gL)$ of $\be(\gL)$ to be dense in $\EE$.

Pull back the uniformity of $\EE$ to $G$ by using the entourages 
\[
    \bigl\{(s,t) \in G\times G : \bigl(\be(s+\gL), \be(t+\gL)\bigr) 
    \in U(F,\varepsilon) \bigr\}.
\]  
This new uniformity on $G$ is compatible with the group structure
(this comes down again to the $G$-invariance of each of the
fundamental entourages) and we have a uniformly continuous mapping of
$G$, equipped with this new topology, into $\EE$. In fact, $\EE$ is
the completion of $G$ under this new uniform topology, and since this
latter is also an Abelian group, $\EE$ becomes a compact Abelian
group.  In more detail, since $G$ is getting its structure by pulling
back the induced structure on $\be(G + \gL)$ in $\EE$ under $t \mapsto
\be(t+\gL)$, we may apply \cite[Ch.~II.3.8, Prop.~~17]{Bou} to see
that the Hausdorff space associated with $G$ under this new uniformity
is homeomorphic to $\be(G + \gL)$ and hence also their completions are
homeomorphic.

So, out of the continuity of the eigenfunctions, we obtain a new
compact Abelian group $\EE$ and a torus parametrization
\begin{equation} \label{XtoE}
   \be \!:\; \XX \longrightarrow \EE \, .
\end{equation} 

By construction, we have:
\begin{prop}\label{separation}
   For each $\lambda\in P(T)$, there exists a unique continuous
   function $g_\lambda$ on $\EE$ with $f_\lambda = g_\lambda \circ
   \be$.  \qed
\end{prop}

Since $(\EE,G,\theta_\EE)$ is pure point with eigenvalues $P(T)$, it
must be conjugate to the $G$-action on $\bS := \widehat{P(T)}$. This,
and in fact a more general statement, is known as the Halmos--von
Neumann representation theorem, compare \cite[Thm.~5.18]{Wal}.  In the
case at hand, we give a short proof.  Explicitly, equip the subgroup
$P(T)$ with the discrete topology, so that its dual group $\bS$ is
compact.  Since $G$ is mapped homomorphically into $\bS$ (namely each
$g \in G$ goes to the evaluation map at $g$ of $P(T)$), it follows
that $\bS$ admits a canonical (and minimal) action of $G$. Fix $x_0\in
\EE$ and normalize the $g_\lambda$ by requiring $g_\lambda(x_0) = 1$,
$\lambda \in P(T)$. Then, $|g_\lambda(x)| = 1$ for all $x \in \EE$.

\begin{prop}  
  The groups $\EE$ and $\bS$ are isomorphic as topological groups by
  the mapping $j:\EE \longrightarrow \bS$ defined by $j(x) \!:\, P(T)
  \longrightarrow U(1), \;\lambda \mapsto g_\lambda(x)$, and thereby
  the dynamical systems $(\EE,G)$ and $(\bS,G)$ are topologically
  conjugate.
\end{prop}
\begin{proof} 
By the normalization condition on $g_\lambda (x_0)$, we infer
\[
    g_\lambda (x) g_\mu (x) \, = \, g_{\lambda \mu} (x)
    \quad \mbox{and} \quad
    g_{\lambda^{-1}} (x) \, = \, \overline{g_\lambda (x)}.
\]
This implies that $j(x)$ is indeed an element of $\bS=\widehat{P(T)}$.
Now, continuity of $j$ follows directly from the continuity of the
$g_\lambda$.  One checks that $j$ is a $G$-map. Injectivity of $j$
follows as the $g_\lambda$, $\lambda \in P(T)$, separate the points of
$\EE$. To show that $j$ is onto, it suffices to show that the dual
$j^*$ of $j$
\[
   j^* \!:\, \widehat{\bS} = P(T) \longrightarrow \widehat{\EE}, \quad
   j^* (\lambda) := \lambda \circ j
\]
is injective (since the image of $j$ is closed and the action of $G$
on $\bS$ is minimal). This can be seen as follows.  Let
$\lambda_1,\lambda_2 \in P(T)$ be given, with
\[
    j^* (\lambda_1) \; = \; j^* (\lambda_2).
\] 
Then, $ \lambda_1 \circ j = \lambda_2\circ j$, i.e., $ j(x)(\lambda_1)
= j(x) (\lambda_2)$ for every $x \in \EE$. As $j(x)(\lambda_1) =
g_{\lambda_1} (x)$ and, similarly, for $ \lambda_2$, we see that this
means $g_{\lambda_1} = g_{\lambda_2}$ which, in turn, implies
$f_{\lambda_1} = f_{\lambda_2}$, and finally $\lambda_1 = \lambda_2$.

Thus, we see that $j$ is indeed a continuous bijection between
compact spaces. Therefore, the inverse of $j$ is continuous as well.

One has to show two more things: that $j$ is compatible with the
group action, and that $j$ is a group homomorphism. Both of them are
more or less straightforward calculations.
\end{proof}

\subsection{Pure point dynamical spectrum together with continuous 
eigenfunctions imply a torus parametrization} \label{special} 
In this section, we specialize the setting of the previous paragraph
by assuming that $\gL$ is Meyer and $(\XX (\gL),G)$ is uniquely
ergodic with pure point dynamical spectrum and continuous
eigenfunctions.  The main objective is to prove Theorem
\ref{EequalA}. Given equation \eqref{XtoE}, it remains to show that
$\EE$ as constructed above is isomorphic to the hull $\bA =\bA(\gL)$
of $\gL$ in the mixed autocorrelation topology. This is done in the
next two paragraphs. Here, we provide some preparation.

\smallskip
As in the proof of Lemma~\ref{lemma.support}, let $C_\mathrm{c}(G)$ denote
 the space of continuous complex-valued functions of compact support
 on $G$ and define, for $c\in C_\mathrm{c}(G)$, the function $\varphi_c
 \!:\, \XX \longrightarrow \CC$ by $\varphi_c(\gL) = \big(c*\delta_
 {\gL}\big)(0)$. Let $g_c := c*\tilde{c} * \gamma^{}_{\gL}$ be the
 corresponding smoothed out autocorrelation of $\gL$, which is a
 continuous function on $G$. From unique ergodicity and Dworkin's
 argument \cite{Dworkin, LMS-1},
\begin{equation} \label{dworkin}
   g_c(t) \; = \; \langle T_{t} \varphi_c, \varphi_c \rangle
\end{equation}
where $t\in G$ and $\langle . , . \rangle$ is the inner product on
$L^2(\XX,\mu)$ whereby it is a Hilbert space. It is not hard to see
that the function $t\mapsto \langle T_{t} f, f\rangle$ is continuous,
bounded and positive definite for any $f\in L^2 (\XX,\mu)$. Thus, by
Bochner's theorem \cite{Loomis}, there exists a unique finite positive
measure $\sigma_f$ on $\Ghat$ with
\[  
    \langle T_{t} f, f\rangle  \; = \; 
    \int_{\Ghat} \, (\widehat{s},t) \dd \sigma_f (\widehat{s})
\]
(see \cite{BL,BL2} as well for a further discussion of this).  It turns
out that, in our context, the spectral measure $\sigma_{\varphi_c}$ can
be explicitly calculated for $c\in C_\mathrm{c} (G)$ in terms of
$\gammahat$. More precisely,
\begin{equation}\label{spectralmeasure}
   \sigma_{\varphi_c} \; = \; |\widehat{c}|^2 \gammahat.
\end{equation}
for any $c\in C_\mathrm{c} (G)$, compare \cite{Martin2,LMS-1,BL}. This
equation links the dynamical spectrum and the diffraction spectrum.

\begin{prop} \label{orthogonal} 
  Let $\lambda$ be an eigenvalue of the uniquely ergodic dynamical
  system\/ $(\XX(\gL),G)$, with associated normalized eigenfunction
  $f_\lambda$.  Let $h\in L^2(\XX(\gL),\mu)$ be arbitrary. Then,
  $\langle h , f_\lambda\rangle =0$ if and only if $\sigma_h
  (\{\lambda\})=0$.
\end{prop}
\begin{proof} By Stone's theorem, compare \cite[Sec.~36D]{Loomis}, there
  exists a projection valued measure
\begin{equation*} 
  E \! : \; \mbox{Borel sets of $\Ghat$} \; \longrightarrow \;
  \mbox{projections on $L^2 (\XX(\gL),\mu)$} 
\end{equation*} 
with $E(B\cap C) = E(B) E(C)$ for $B,C\subset \Ghat$ measurable  and 
\begin{equation} \label{spectralfamily}
  \langle T_t f, g \rangle \; = \; 
  \int_{\Ghat} (\hat{s},t) \dd \sigma_{f,g} (\hat{s}) 
\end{equation} 
where the measure $\sigma^{}_{f,g}$ on $\Ghat$ is defined by $\sigma^{}_{f,g}
(B) := \langle E (B) f, g \rangle$.  In particular, we have $\sigma^{}_{f,f} =
\sigma_f$.  

\smallskip
From $E(B) E(\{\lambda\}) = E(B\cap \{\lambda\})$, we infer that
$\sigma^{}_{E(\{\lambda\})f,g}$ is concentrated on $\{\lambda\}$ for
arbitrary $f,g\in L^2 (\XX(\gL),\mu)$.  This easily yields that $T_t
E(\{\lambda\}) f = (\lambda,t) E(\{\lambda\}) f$ for any $f\in L^2
(\XX(\gL),\mu)$.

\smallskip 
Conversely, if $f$ is an eigenfunction to $\lambda$, we infer from the
validity of the equation $(\lambda,t) \langle f,g\rangle = \langle T_t
f,g\rangle = \int_{\Ghat} (\hat{s},t) \dd \sigma_{f,g} (\hat{s}) $,
for all $t\in G$, that $\sigma_{f,g}$ is concentrated on $\lambda$. This
easily gives $E(\{\lambda\}) f = f$ for any eigenfunction to the
eigenvalue $\lambda$.

\smallskip
Put together, this means that $E(\{\lambda\})$ is the orthogonal
projection onto the eigenspace for the eigenvalue $\lambda$. This
eigenspace is one-dimensional by ergodicity. Thus, $ E(\{\lambda\}) h =
\langle h, f_\lambda\rangle f_\lambda$ and we infer
\[ 
   |\langle h , f_\lambda\rangle|^2 \, = \,
    \|E(\{\lambda\}) h\|^2 \, = \, 
    \langle E(\{\lambda\}) h, E(\{\lambda\}) h\rangle 
    \, = \, \sigma_h (\{\lambda\}).
\]
Now, the statement of the proposition is immediate.
\end{proof}

We are now ready to prove the isomorphism between $\EE$ and
$\bA=\bA(\gL)$.  Both spaces in question, $\EE$ and $\bA$, are
obtained by completion of $G$ in uniform topologies for which a
fundamental system of $G$-\emph{invariant} entourages exist. For this
reason, it is actually sufficient to show that the identity mapping
from $G$ to itself is bi-continuous at $0$ when these two topologies
are put on two sides.

\subsection{Continuity of $\bA \longrightarrow \EE$}

If $\gL$ is a Meyer set, we know that $\varDelta=\gL - \gL$ is
uniformly discrete. Consequently, there is a compact neighbourhood $K$
of $0$ in $G$ so that, for all $t \in \varDelta$, $(t+ K)\cap
\varDelta = \{t\}$.

Let $\{x_i\}$ be a net in $G$ which converges to $0$ in the
autocorrelation topology. Then, there are elements $v_i \in G$
converging to $0$ in the original topology of $G$ so that $d(v_i + x_i
+ \gL, \gL)$ converges to $0$. Here, $d$ is the pseudo-metric defined
in \eqref{metric}, which, by \eqref{metricvseta}, satisfies
\[ 
   d(s+\gL, t +\gL) \; = \;
   \lim_{n \to \infty}\frac{ \mathrm{card} 
   \left( \big((s+\gL) \, \symdiff \, (t + \gL)\big) 
   \cap A_n \right)}{\theta_G(A_n)} \; = \; 
   2 \big(\eta(0) - \eta(s-t)\big) \, ,
\]
so, for $y_i:= v_i+x_i$,
\begin{equation}\label{dtozero}
   d(y_i +\gL, \gL) \; = \; 2 \big(\eta(0) - \eta(y_i)\big)
   \;\rightarrow\; 0 \, .
\end{equation}
This convergence of the $\{y_i +\gL\}$ to $\gL$ shows that $y_i \in
\varDelta$ for all sufficiently large $i$. If we show that $\{y_i\}$
converges to $0$ in the topology of $\EE$, this will also give
convergence of the original net $\{ x_i \}$, since the topology of
$\XX (\gL)$, and hence $\EE$, is defined so that shifts by small
elements of $G$ are small.

Now, $\gamma^{}_{\gL} = \sum_{t \in \varDelta} \eta(t)\delta_t$, 
Let $c \in C_\mathrm{c} (G)$ with $\supp(c*\tilde c)\subset K$.
By our choice of $K$, $0$ is then the only element of
  $\varDelta$ in $\supp (c*\tilde c)$. Thus,
\[
   g_c(y_i) \, = \, \big( c*\tilde c*\gamma^{}_{_\gL} \big)(y_i) 
   \, = \, \sum_{t \in \varDelta} \eta(t)\,
   \big(c*\tilde{c}\big) (y_i - t) 
   \, = \, \eta(y_i)\, \big(c*\tilde{c}\big)(0) \, .
\] 
By \eqref{dtozero}, this implies $\lim_i g_c(y_i) = g_c(0)$. By
\eqref{dworkin}, this means that
\begin{equation} \label{converge}
    \langle T_{y_i} \varphi_c , \varphi_c\rangle 
    \;\longrightarrow\; \langle \varphi_c, \varphi_c \rangle.
\end{equation}  

As we have pure point spectrum with the set of eigenvalues $P(T)$ and
corresponding normalized eigenfunctions $f_\lambda$, $\lambda\in
P(T)$, we can write $\varphi_c$ as a Fourier series, $\varphi_c =
\sum_{\lambda\in P(T)} a_\lambda f_\lambda$, where the $a_\lambda$
(which depend on $c$) are complex numbers.

Then, using \eqref{converge}, we find 
$ \sum \lambda(y_i) \,| a_\lambda |^2  \| f_\lambda \|^2 \longrightarrow
\sum |a_\lambda |^2 \| f_\lambda \|^2 \, $ which results in 
\[
     \sum |a_\lambda|^2 \big(\lambda(y_i) -1 \big) 
     \;\longrightarrow\; 0 \, , 
\]
by normalization of the eigenfunctions. Taking complex conjugates
then yields
\[
    \sum |a_\lambda|^2(\overline{\lambda(y_i)} -1) 
    \;\longrightarrow\;  0 \, . 
\]
As $\lambda$ takes values in $U(1)$, we have
 \[ 
    |\lambda(y_i) -1|^2 \; = \;  (1 -
     \lambda(y_i)) + (1 - \overline{\lambda(y_i)}) 
\]
and  we obtain
\[ 
   \sum |a_\lambda|^2 |\lambda(y_i) -1|^2  \; = \; 
   \sum |a_\lambda|^2 \big( 1 - \lambda(y_i) +
   1 - \overline{\lambda(y_i)} \big) 
   \; \longrightarrow \; 0 \, .  
\]
Thus, $\{\lambda(y_i)\} \to 1$, whenever $a_\lambda\neq 0$. Now,
$a_\lambda\neq 0$ means that $\varphi_c$ is not orthogonal to
$f_\lambda$.  By Proposition~\ref{orthogonal}, this is equivalent to
to $\lambda \in P(\sigma^{}_{\varphi_c})$ (see \eqref{defOfP} for the
definition of $P( \cdot)$).  Thus, we have $\{\lambda(y_i)\} \to 1$
for all $\lambda \in P(\sigma_{\varphi_c})$ and for all $c\in C_{\rm
  c}(G)$.  As $c\in C_\mathrm{c}(G)$ is arbitrary,
\eqref{spectralmeasure} then shows that $\{\lambda(y_i)\} \to 1$ for
all $\lambda \in P(\widehat\gamma)$, and then for all $\lambda \in
\langle P(\widehat\gamma) \rangle$. As $ P(\widehat\gamma) = P(T)$ by
Lemma~\ref{lemma.purepoint}, this means, for all $\lambda \in P(T)$,
\begin{equation*}
    f_\lambda(y_i +\gL) \; = \; \lambda(y_i)f_\lambda(\gL) 
   \; \longrightarrow\; f_\lambda(\gL)
\end{equation*}
and this is precisely the meaning of convergence of $\{y_i\}$ to $0$
in the $\EE(\gL)$-topology. This concludes the continuity
argument in the first direction.

\subsection{Continuity of $\EE \longrightarrow \bA$}

Let $\lambda \in P(T)$, so $f_\lambda(x+\gL) = \lambda(x)
f_\lambda(\gL)$ for all $x\in G$.  The continuity of $f_\lambda$ then
shows that $|f_\lambda | $ is a non-zero constant function on $\XX
(\gL)$.  Let $\{x_i\} \to 0$ in the $\EE$-topology on $G$. Then,
$\{f_\lambda(x_i+\gL)\} \to f_\lambda(\gL)$ shows that
$\{\lambda(x_i)\} \to 1$ for all $\lambda \in P(T)$.

Let $c \in C_\mathrm{c}(G)$ be chosen so that $0 \le c(x) \le 1$ for all
$x\in G$, with $c(x) =1 \Leftrightarrow x = 0$. Moreover, let $c$ be
so that $\nu:= c*\tilde c$ satisfies $(\supp(\nu) - \supp(\nu)) \cap
(\varDelta - \varDelta) = \{0\}$, which rests upon the Meyer
property. Then, $ \|c\|_2^2 = \nu(0) > \nu(x) \ge 0$ for all $x \in
G\backslash\{0\}$ and $\supp(\nu) \cap \varDelta = \{0\}$.

Let $\varphi_c = \sum_{\lambda \in P(T)} a_\lambda f_\lambda$ be the
Fourier expansion of $\varphi_c$.  Choose $\varepsilon >0$ and find a
finite set $F \subset P(T)$ so that
\[
   \| \varphi_c - \sum_{\lambda \in F} a_\lambda f_\lambda\|_2 
  \; < \; \varepsilon \, .
\] 
Choose $N$ in the index set of $\{x_i\}$ so that, for all $i
\succcurlyeq N$ and all $\lambda \in F$,
\[
   |\lambda(x_i) -1 | \; < \; \frac{\varepsilon}{1 + 
       \sum_{\lambda\in F} |a_\lambda|} \, .
\] 
Then,
\begin{eqnarray*}
   \|T_{x_i} \varphi_c - \varphi_c\|_2 & < & 
   \big\|T_{x_i} \varphi_c - T_{x_i} \sum_{\lambda \in F} 
                 a_\lambda f_\lambda\big\|_2 \\ 
  &+& \big\|\sum_{\lambda \in F} \lambda(x_i)\, a_\lambda  f_\lambda 
            - \sum_{\lambda \in F} a_\lambda f_\lambda\big\|_2 + 
      \big\|\sum_{\lambda \in F} a_\lambda f_\lambda - \varphi_c\big\|_2\\
  & < & \varepsilon + \sum_{\lambda \in F} |\lambda(x_i) -1|
    \, |a_\lambda| \,  \|f_\lambda\|_2 + \varepsilon  
    \; < \; 3\varepsilon \, ,
\end{eqnarray*}
since the $T_x$ are unitary and $\|f_\lambda\|_2 = 1$.

Thus, $\{T_{x_i} \varphi_c\} \to \varphi_c$ and hence $\{\langle
T_{x_i} \varphi_c, \varphi_c\rangle \} \to \langle \varphi_c,
\varphi_c \rangle$ and, using (\ref{dworkin}) again,
\begin{equation} \label{conv}
   \bigl\{ g_c(x_i) = \sum_{t\in\varDelta} \nu(x_i -t) \eta (t) \bigr\} 
   \;\longrightarrow\; g_c(0) \, .
\end{equation}
{}For each $x_i$, there is at most one $t_i \in \varDelta$ with
$x_i - t_i \in \supp(\nu)$. Thus,
\begin{equation*}
    g_c(x_i) \; = \; \begin{cases}
      \nu(x_i - t_i)\eta(t_i), & \text{if $t_i$ exists,}\\
         0, &\text{otherwise.}
\end{cases}
\end{equation*}
Moreover, $g_c (0) = \nu(0)\eta(0)$.

Now, \eqref{conv} implies that $\nu(x_i - t_i)\eta(t_i) \to \nu(0)
\eta(0) \ne 0$ (so, in particular, the $t_i$ must exist
eventually). Since $0 \le \nu(x_i - t_i) \le \nu(0)$ and $0 \le
\eta(t_i) \le \eta(0)$, we get $\{\eta(t_i)\} \to \eta(0)$ and
$\{\nu(x_i - t_i)\} \to \nu(0)$.  By the choice of $c$, $\{v_i:= t_i -
x_i\} \to 0$.

Now, $\{x_i\}$ converges to $0$ in the $\bA$-topology since $\{v_i +
x_i\} = \{t_i\}$, the $\{v_i\} \to 0$ in the original topology of $G$,
and $\{ d(t_i +\gL, \gL) = 2(\eta(0) -\eta(t_i))\} \to 0$.

This implies continuity in the other direction and completes the proof
of Theorem ~\ref{EequalA}. \qed

\subsection {Proof of Theorem~$\ref{main4}$, (b) $\Longrightarrow$ (a)}
This is immediate from  Theorem~\ref{EequalA}.  \qed

\section{The proof of the sufficiency direction of Theorem~\ref{main1}}
\label{Proof_of_main1}

\begin{proof}
  The hypotheses of Theorem~\ref{main1}, in the direction
    of sufficiency, include those of Theorem~\ref{main4}. Thus, we
  have a torus parametrization $\XX(\gL) \longrightarrow \bA(\gL)$
  coming from the mapping $\XX(\gL) \longrightarrow \EE$.  This
  provides us with a cut and project scheme according to
  Theorem~\ref{prop.torusparametrization}.

By assumption, our torus parametrization is one-to-one almost
everywhere. Thus, Theorem~\ref{intermediate} and
Theorem~\ref{prop.measureboundary} apply.  Therefore, we obtain a
window by Theorem~\ref{intermediate}, whose boundary has Haar measure
$0$ in $H$ by Theorem~\ref{prop.measureboundary}.

Repetitivity of $\gL$ is equivalent to the minimality of $\XX(\gL)$,
which is assumed. Thus, $\XX(\gL)$ is associated with a regular model
set according to Theorem~\ref{intermediate} as required.
\end{proof}

\section{The proof of necessity in Theorem~\ref{main1}} \label{toruscp}

Up to now, the direction of investigation has been from dynamical
systems and torus parametrizations to cut and project schemes and
model sets. In this section, we derive results going in the other
direction and prove that the conditions of Theorem~\ref{main1} are
necessary.

\subsection{Irredundancy}
Assume that we are given an IMS $\gL$ with respect to the CPS
$(G,H,\LL)$. We wish to construct a torus parametrization for the
local hull of $\gL$. As we have seen in Proposition~\ref{beta-prop},
the curious property of irredundancy is crucial to the existence of
such a map. What happens if we have a CPS together with an IMS for
which the window fails irredundancy?  Can we modify the CPS (and
thereby the torus) to get the irredundancy? The answer is yes. But 
it is here that the notion of an \emph{inter} model set becomes
important.

\begin{lemma} \label{quotient} 
  Let $(G,H,\LL)$ be a CPS. Let $W\subset H$ be a non-empty compact
  set with $\overline{W^\circ}=W$ and\/ $\theta_H (\partial W)=0$. Let
  $\gL\subset G$ with $\oplam (W^\circ) \subset \gL \subset \oplam
  (W)$ be arbitrary. Then, there exists a CPS\/ $(G,H',\LL')$ and\/
  $W'\subset H'$ compact, non-empty and irredundant, with
  $W'=\overline{{W'}^\circ}$ and $\theta_{H' }(\partial W') =0$, such
  that $\oplam ({W'}^\circ) \subset \gL \subset \oplam (W')$, i.e.,
  that $\gL$ is a regular IMS for the new CPS\/ $(G,H',\LL')$ with
  window $W'$.
\end{lemma}

\begin{remark} 
  The proof of the lemma relies on factoring out the stabilizer of
  $W$. It is crucial to note that sets of the form $\oplam (W)$ may
  \emph{not} be representable as sets of the form $\oplam (W')$ in
  the emerging ``quotient'' scheme.  Rather, sets lying between
  $\oplam ({W}^\circ)$ and $\oplam (W)$ can be exhibited as sets lying
  between $\oplam({W'}^\circ)$ and $\oplam(W')$. We refer the reader
  to \cite{LM} for further discussion.
\end{remark}

\begin{proof}
Let $(G,H,\LL)$ be the given CPS, with the usual conditions on
the projections $\pi^{}_{1}$ and $\pi^{}_{2}$, and with the compact
regular window $W=\overline{W^\circ}\neq \varnothing$.  Let 
\[
   I \; := \; \mathrm{stab}^{}_H (W) \; = \; \{ t \in H : t+W=W \} 
\]
be the stabilizer of $W$, which is a subgroup of $H$. Clearly, $I$ is
closed, and $I\subset W-W$ implies that $I$ is compact. Observe that
we also have $W^\circ + I = W^\circ$.

Define the factor group $H' = H/I$ and let $\rho \! : \; H
\longrightarrow H'$ be the natural map. Moreover, define
\[
   \LL' \, := \, \{ (x, \rho(x^\star)) : x\in L\}
   \; \subset \; G\times H'
\]
together with a mapping $\LL\longrightarrow \LL'$ defined
by $(x,x^\star)\mapsto (x,\rho(x^\star))$. This is a group
homomorphism, and surjective.  Since $(x,x^\star) \mapsto (0,0)\in
\LL'$ is only possible for $x=0$, we see that also
$\rho(x^\star)=0$ in this case.  Consequently, the kernel of this
homomorphism is $\{(0,0)\}$, and $\LL \simeq \LL'$.

Consider the diagram
\[ \begin{CD}
   G\times H @>{\text{id}\times\rho}>> G\times H' \\
   @V\text{nat}VV @VV\text{nat}V  \\ 
   (G\times H)/\LL @>>> (G\times H')/\LL'
  \end{CD}
\]
where the horizontal arrow in the lower line exists in an obvious way,
because $(x,x^\star)\in \LL$ is mapped to $(x,\rho(x^\star))\in\LL'$.
Since $\LL$ is a closed subgroup of $G\times H$, quotient is
Hausdorff.  Moreover, the natural mapping $G\times H \longrightarrow
(G\times H)/\LL$ is an open map, and we get that $(G\times H)/\LL
\longrightarrow (G\times H')/\LL'$ is continuous.

Consequently, $(G\times H')/\LL'$ is compact, whence
$\LL'$ is co-compact. Consider now a compact neighbourhood
$U\times \rho(V)$ of $0$ in $G\times H'$, where $V$ is a compact
neighbourhood of $0$ in $H$ and $U$ is a compact neighbourhood of $0$
in $G$. Then,
\[
   \LL' \cap (U\times \rho(V)) \; = \;
   \{(x,\rho(x^\star)) : x\in U, \, \rho(x^\star)\in V\}
   \; = \; \{(x,x^\star) : x\in U, \, x^\star\in V+I\} \, .
\]
Since $U\times(V+I)$ is compact and $\LL$ is a lattice, the set
$\{(x,x^\star) : x\in U, \, x^\star \in V+I\}$ is \emph{finite}, and
contains $(0,0^\star) = (0,0)$. Consequently, $\LL'\cap
\big(U\times\rho(V)\big)$ is finite, too, 
with $(0,0)\in \LL'\cap \big(U\times\rho(V)\big)$. Consequently,
$(0,0)$ is isolated and $\LL'$ is (uniformly) discrete.

This shows that $(G,H',\LL')$ is another CPS, with all the properties
required, for which we now need a window. To this end, define $W' =
\rho(W)$. We note that $W'$ is compact since $W$ is.
  Moreover, also $W$ is the complete preimage of $\rho(W)$ since $W=
  I+W$.  As $\rho$ is an open map, $\rho(W^\circ)$ is open and we have
  ${W'}^\circ = \rho(W^\circ)$ and $W' = \overline{{W'}^\circ} \neq
  \varnothing$.

Note that $\rho \! : \, H \longrightarrow H'$ induces a mapping from
the Haar measure $\theta^{}_{H}$ to a Haar measure $\theta^{}_{H'}$,
with
\[
   \theta^{}_{H'} (\partial W') \; = \; 
   \theta^{}_{H'} (W' \setminus {W'}^\circ ) \; = \;
   \theta^{}_{H} \big( (W+I) \setminus (W^\circ +I)\big) \; = \;
   \theta^{}_{H} (W\setminus W^\circ) \; = \; 0 \, .
\]
Finally, if $\oplam(W^\circ) \subset \gL \subset \oplam(W)$ 
for some $\gL\subset G$ in the original CPS, then
\[
   \oplam({W'}^\circ) \; \subset \; \gL
   \; \subset \; \oplam(W')
\]
in the new CPS. This is easy to check because $x\in\oplam({W'}^\circ)
\iff \rho(x^\star)\in {W'}^\circ \iff x^\star \in W^\circ$ and
similarly for the remaining details.  This completes the argument.
\end{proof}

\subsection{Torus parametrizations for model sets}

Assume that we are given an IMS $\gL$ with respect to the CPS
$(G,H,\LL)$ with a window $W$. According to the previous paragraph,
we may assume that the CPS is irredundant, and then, by
Proposition~\ref{beta-prop}, that we have a torus parametrization of
the local hull into the compact group $\TT$ of this scheme. We shall
now convert this to a torus map into $\bA(\gL)$.

\begin{lemma} \label{lemma.aux}  
   Let\/ $(G,H,\LL)$ be a CPS and let\/ $W\subset H$ be a compact,
   non-empty and irredundant window, with $W=\overline{W^\circ}$ and\/
   $\theta_H (\partial W) =0$. Then, the following holds.
\begin{itemize}
\item[\textrm{(a)}] $\theta_H ((W - h_n) \symdiff W)
   \longrightarrow 0$, whenever $\{h_n \}$ is a net in $H$ with 
   $h_n \to 0\in H$.  
\item[\textrm{(b)}] If $h\in H$ satisfies $\theta_H ((W-h)
   \symdiff W) = 0$, then $h=0$.
\end{itemize}
\end{lemma}
\begin{proof} 
(a) Denote the canonical representation of $H$ on $L^1 (H,\theta_H)$
by $\tau^H$, i.e., $\tau^H_h f (x) = f( -h +x)$. It is well known that
this representation is strongly continuous \cite{Rud}. This means that
$\tau^H_h f \longrightarrow f$ for $h \to 0$ and all $f\in L^1
(H,\theta_H)$. Now, let $1_W$ be the characteristic function of the
compact set $W\subset H$. Then, $1_W$ belongs to $L^1 (H,\theta_H)$
and therefore
\[ 
    \theta_H ((W - h) \symdiff W) \; = \;
     \| \tau^H_h 1_W - 1_W\|_{L^1} 
     \;\longrightarrow\; 0 \, , \quad \mbox{as }
     h \to 0 \, .
\] 
This proves (a).

\smallskip
(b) Let $h\in H$ be given with $\theta_H ((W -h) \symdiff W) = 0$.  As
the window is translationally fixed only by $0\in H$, it suffices to
show $W-h =W$, i.e., $(W -h ) \symdiff W =\varnothing$. Assume the
contrary. Then, $(W -h ) \symdiff W$ actually contains an open set
because $W = \overline{W^\circ}$. This implies $\theta_H ( (W -h )
\symdiff W ) >0$ and we arrive at a contradiction.
\end{proof}

\begin{prop} \label{continuityofj}   
  Let $(G,H,\LL)$ be a CPS with associated torus\/ $\TT$. Let\/
  $W\subset H$ be a compact, non-empty and irredundant window, with
  $W=\overline{W^\circ}$ and $\theta_H (\partial W) =0$. Let
  $\gL\subset G$ be an arbitrary IMS, i.e., $\oplam (W^\circ) \subset
  \gL \subset \oplam (W)$. Then, selecting a  
  neighbourhood $U$ so that $\gL\in \CalD_U$, the mapping $j \! : \, \TT
  \longrightarrow \CalDq_U$, $(x,h) + \LL \mapsto \beta(x
  +\oplam(W-h))$, is continuous and injective. In particular, $j$
  gives an isomorphism between $\TT$ and $\bA (\gL)$.
\end{prop}

\begin{proof} 
We first show continuity of $j$. Let $\{\xi_\iota\}$ be a net in $\TT$
with $\xi_\iota\longrightarrow \xi\in \TT$.  Then, without loss of
generality, we may assume that $\xi = (x,h) + \LL$, $\xi_\iota =
(x_\iota,h_\iota) + \LL$ and $x_\iota \longrightarrow x$, $h_\iota
\longrightarrow h$. As the topology of $\CalDq_U$ allows for small
translations, it suffices to show that $\beta (\oplam(W-
h_\iota))\longrightarrow \beta (\oplam (W-h)).$ By part (a) of
Lemma~\ref{lemma.aux}, we infer
\begin{equation}\label{convergencetozero}
   \theta_H ((W-h_\iota) \symdiff (W-h)) \;\longrightarrow\; 0 \, .
\end{equation}
Therefore,
\begin{eqnarray*}
   d \big(\beta (\oplam(W- h_\iota) ) , \beta (\oplam(W- h) ) \big) 
  &=& \dens \big(\oplam(W- h_\iota) \symdiff \oplam(W- h) \big)\\ [1mm]
  &=& \dens \big( \oplam( (W- h_\iota) \symdiff (W- h) ) \big) \\ 
  &\stackrel{\mathrm{by\; Thm.~\ref{uniformdistribution}}}{=}&
    \mathrm{dens} (\LL)\, \theta_H ((W-h_\iota)\symdiff (W-h))\\
  &\xrightarrow{\mathrm{by\; \eqref{convergencetozero}}} & 0.
\end{eqnarray*}
This proves the continuity statement.

\smallskip
We now show injectivity. To do so, let $(x,h)$ and $(x',h')$ in
$G\times H$ be given with \newline $\beta (x + \oplam (W -h)) = \beta
(x' + \oplam (W -h'))$. This implies
\[ 
   \beta ( \oplam( z^\star + k +  W' )) \; = \; \beta (\oplam(W')) \, ,
\] 
where $ z= x - x'$, $k = h' -h$ and $W' = - h' + W$. By
Proposition\ref{beta-prop} and
Theorem~\ref{uniformdistribution}, we then obtain
\begin{eqnarray*}
   \mathrm{dens}(\LL)\,
   \theta_H \big((z^\star + k +  W')\symdiff \oplam(W') \big)
   & = & \dens \big(\oplam (z^\star + k +  W' \symdiff \oplam(W') \big) \\
   & = & d \big(\beta (\oplam( z^\star + k +  W')) , \beta (\oplam(W')) \big) 
   \; = \; 0.
\end{eqnarray*}
By part (b) of Lemma~\ref{lemma.aux}, this gives $0 = k +z^\star$
or, put differently, $(x,h)+ \LL = (x',h') +
\LL$. This proves injectivity.

\smallskip
The inverse of a continuous injective map on a
compact space is continuous as well. 

\smallskip
So far, we know that $j$ is a continuous bijective map from $\TT$ onto
$j(\TT)$. By continuity of $j$ and minimality of the action of $G$ on
$\TT$, $j(\TT)$ is just $\bA (\oplam (W))$. By uniform distribution,
$\beta(\oplam (W)) = \beta(\gL)$ and $\bA (\oplam (W))= \bA (\gL)$
follows. This proves the last statement.
\end{proof}

\begin{remark}
It is known from \cite{MS} that $ \TT$ and $\bA (\gL)$ are isomorphic
for regular model sets, and our proof is in some sense a
variant of the proof in \cite{MS}. However, our result here is more
explicit in that $\beta$ is shown to establish this isomorphism, and
this will allow us to clarify the relationship between $\beta$ and
$\ba$.
\end{remark}

\begin{theorem}\label{main3} 
  Let\/ $(G,H,\LL)$ be a CPS, and let a non-empty window $W \subset H$
  with\/ $W = \overline{W^\circ}$ and $\theta_H (\partial W)=0$ be
  given. If\/ $\gL \subset G$ satisfies\/ $t+ \oplam(W^\circ) \subset
  \gL \subset t+ \oplam(W)$ for some\/ $t\in G$, then the canonical
  mapping\/ $\beta \!:\; \XX (\gL) \longrightarrow \bA (\gL)$ is
  continuous and\/ one-to-one almost everywhere.
\end{theorem}
\begin{proof}
  By Lemma~\ref{quotient}, we may assume, without loss of generality,
  that the cut and project scheme is irredundant. Assume furthermore
  that $t=0$.  Proposition~\ref{beta-prop} then gives a torus
  parametrization $ \bt \! : \; \XX(\gL) \longrightarrow \TT $ with
\[ 
   \bt (\gL')  =  (x,h) + \LL 
   \; \Longleftrightarrow \;
    x + \oplam(-h + W^\circ) \subset \gL'
\subset x + \oplam (-h + W).
\]
By $\theta_H (\partial W)=0$ and Theorem \ref{uniformdistribution},
$\beta (\gL') = \beta( x +\oplam(-h + W))$ whenever $x + \oplam(-h +
W^\circ) \subset \gL' \subset x + \oplam (-h + W)$.  Thus,
Proposition~\ref{continuityofj} shows that $\beta = j\circ \bt$ with a
continuous injective $j$. Thus, $\beta$ is a continuous torus
parametrization.  It remains to show that it is one-to-one almost
everywhere. But this follows from Theorem~\ref{prop.measureboundary}.
\end{proof}

\begin{remark} 
  Proposition~\ref{continuityofj} and the considerations in the proof
  of the previous theorem effectively show that the maps $\ba \!: \,
  \XX (\gL) \longrightarrow \bA$ and $\bt \! : \, \XX
  (\gL)\longrightarrow \TT$, yielding a description of $\gL$ as a
  regular model set in the previous sections, can be identified with
  the canonical map $\beta$.
\end{remark}
\smallskip

\subsection{The end of Theorem~\ref{main1}}
\label{proofs}

\begin{proof} 
  Let $\gL$ be a regular model set. By Lemma~\ref{quotient}, we may
  assume that its CPS is irredundant and that we are in the situation
  of Theorem~\ref{main3}. This provides us with an almost everywhere
  one-to-one continuous mapping of $\XX(\gL)$ onto $\bA(\gL)$.
  Theorem~\ref{main2} shows that $\XX(\gL)$ is uniquely ergodic, and,
  since $\gL$ is repetitive, $\XX(\gL)$ is also minimal. By
  Theorem~\ref{main4}, we obtain a pure point dynamical spectrum with
  continuous eigenfunctions that separate almost all points of
  $\XX(\gL)$.
\end{proof}

\section{The crystallographic case} \label{crystallographic}

The aim of this section is to give a proof of the following result
on the characterization of fully periodic Delone sets.

\begin{theorem}\label{main5}
   Let\/ $G$ be an LCA group and $\gL$ a uniformly discrete subset 
   of $G$. Then, the following assertions are equivalent.
\begin{itemize}
\item[(i)] $\gL$ is  crystallographic.
\item[(ii)] $\gL$ is Meyer and the map $\beta :\; 
            \XX (\gL)\longrightarrow \bA(\gL)$ is 
            continuous and injective. 
\item[(iii)] All of the  following conditions hold:
\begin{itemize}
\item[(1)] All elements of\/ $\XX(\gL)$ are Meyer sets. 
\item[(2)] $(\XX(\gL),G)$ is uniquely ergodic.
\item[(3)] $(\XX(\gL),G)$ has pure point dynamical spectrum 
           with continuous eigenfunctions.
\item[(4)] The eigenfunctions separate all points of\/ $\XX (\gL)$. 
\end{itemize}
\end{itemize}
In this case, $(\XX(\gL),G)$ is also minimal, hence strictly ergodic.  
\end{theorem}

Recall that $\gL$ is called \emph{crystallographic} 
(or fully periodic) if its periods 
\[
   \mathrm{per}(\gL) \; := \; \{t\in G : t + \gL = \gL\}
\] 
form a lattice, i.e., a co-compact discrete subgroup of $G$.

\smallskip
We start by proving the equivalence of claims (i) and (ii)
of Theorem~\ref{main5}. 

\begin{lemma} 
   The Delone set $\gL \subset G$ is crystallographic if and only
   if $\gL$ is Meyer and the mapping $\beta \!:\,
   \XX(\gL)\longrightarrow \bA(\gL)$ is continuous and
   injective.
\end{lemma}
\begin{proof} 
The `only if' part is easy: If the Delone set $\gL$ is
crystallographic, with lattice of periods
$P=\mathrm{per}(\gL)$, it is of the form $\gL = F + P$,
where $F$ must be a finite set due to the uniform discreteness of
$\gL$. From $P-P=P$, we have $\gL -
\gL = (F-F) + P$.  This is still uniformly discrete, because
$F-F$ is still finite, and $\gL$ is Meyer.

Let $\TT=G/\mathrm{per}(\gL)$ be the compact quotient of $G$ by the
set of periods.  Then, there is a natural isomorphism
$\TT\longrightarrow \XX (\gL)$.  This easily yields the statement
about $\beta$.

\smallskip
We now prove the `if' statement: As $\beta$ is continuous,
$\bA(\gL)$ is compact and we have pure point diffraction. In
particular, the $\varepsilon$-almost periods $P_\varepsilon$ are
relatively dense.

Obviously, $\mathrm{per} (\gL)$ is a subgroup of $G$. As
  $\gL$ is a Meyer set, $\mathrm{per} (\gL)$, being a subset of $\gL -
  \gL$, is uniformly discrete. Therefore, it suffices to show that
the set of periods is relatively dense (which implies that the
quotient $G/\mathrm{per}(\gL)$ is compact).  We shall show that the
set of periods contains the set of $P_\varepsilon$ of
$\varepsilon$-almost periods for a suitable $\varepsilon >0$. As
$P_\varepsilon$ is relatively dense, the desired statement follows.
Here are the details:

As $\bA (\gL)$ is compact and $\beta$ is continuous 
and injective by assumption, it  has a continuous inverse
\[
   \alpha \!:\, \bA (\gL) \longrightarrow \XX (\gL).
\] 
As $\bA (\gL)$ is compact, $\alpha$ is even uniformly continuous.
Thus, for every compact $C\subset G$ and open $V\subset G$, there
exists an $\varepsilon>0$ such that
\begin{equation}\label{close} 
  (\alpha (\xi) + t,\alpha(\xi)) =  
   (\alpha (\xi + t),\alpha (\xi)) \in U_\mathrm{LT} (C,V)
\end{equation}
for all $\xi \in \bA$ and all $t\in P_\varepsilon$, where use
the addition of $t$ for the translation action on both spaces
for simplicity. Here, of course,
\[
   U_\mathrm{LT} (C,V) \; = \;
  \{(\gG,\gG') : (\gG + v)\cap C = \gG'\cap C \mbox{
  for a suitable $v \in V$} \}.
\]

Now, choose an open neighbourhood $V$ of $0$ in $G$ according to
Fact~\ref{snapin}, meaning that we have $V \cap ((\gL - \gL) + (\gL - \gL))=
\{0\}$, and a compact set $C\subset G$ such that, for all
$t\in G$, $(t + C) \cap \gG \neq \varnothing$.

Choose $\varepsilon>0$ so that \eqref{close} holds. As $\alpha$ is
onto, we infer that, for every $\gG\in \XX (\gL)$ and every $t\in
P_\varepsilon$,
\[ 
(\gG + t, \gG) \in  U_\mathrm{LT} (C,V).
\]
Fact~\ref{snapin} then implies
\[ 
    (\gG + t) \cap C = \gG \cap C
\] 
for every $\gG \in \XX (\gL)$ and $t\in P_\varepsilon$. As
$\gG$ is arbitrary, we can replace $\gG$ by $s +\gG$
with $s\in G$ arbitrary. Thus, we end up with
\[ 
   (\gG + s + t) \cap C = (\gG + s) \cap C
\] 
for every $s\in G$. As $C\neq\varnothing$, this shows that every $t\in
P_\varepsilon$ is a period of $\gG$.
\end{proof}

We now show equivalence of (ii) and (iii). 
\begin{lemma}
   A set $\gL\subset G$ is Meyer and $\beta \!:\, \XX (\gL)
   \longrightarrow \bA (\gL)$ is continuous and injective if and only
   if the following\/ $4$ conditions hold.
\begin{itemize}
\item[(1)] All elements of\/ $\XX(\gL)$ are Meyer sets;
\item[(2)]  $(\XX(\gL),G)$ is uniquely ergodic;
\item[(3)] $(\XX(\gL),G)$ has pure point dynamical spectrum with 
           continuous eigenfunctions;
\item[(4)] The eigenfunctions separate all points of\/ $\XX (\gL)$. 
\end{itemize}
\end{lemma}
\begin{proof} We first show the ``only if'' part: 
The validity of (1) is clear as $\gL$ is Meyer. By assumption on
$\beta$, $\XX (\gL)$ is isomorphic to the group $\bA (\gL)$.  Now,
$(\bA (\gL),G)$ is uniquely ergodic (as the action of $G$ is minimal
and $\bA (\gL)$ is a group) and it has pure point dynamical spectrum
with continuous eigenfunctions given by the characters. As, by
assumption on $\beta$, the dynamical system $(\XX (\gL),G) $ is
topologically conjugate to $(\bA (\gL),G)$, the assertions (2), (3)
and (4) follow.

\smallskip
We now prove the ``if'' part: We can apply Theorem~\ref{main4} to
obtain a continuous mapping $\ba \!:\, \XX (\gL)\longrightarrow \bA (\gL)$
as we have pure point dynamical spectrum with continuous
eigenfunctions. Moreover, again by Theorem~\ref{main4}, the map is
actually injective because the eigenfunctions separate \emph{all} points.

\smallskip
It remains to show that $\ba$ agrees with $\beta$. We first note that
$(\XX (\gL),G)$ is minimal, because $(\bA (\gL), G)$ is minimal and
$\ba$ is injective.  By injectivity of $\ba$ and
Theorem~\ref{intermediate}, $(\XX (\gL),G)$ is associated with a
repetitive model set. In fact, the model set is regular by
Theorem~\ref{prop.measureboundary}.  Thus, by Theorem~\ref{main3},
the map $\beta$ is continuous.  Consequently, $\beta$ and $\ba$ are
continuous $G$-maps from $\XX (\gL)$ into $\CalDq_{U}$, for suitable
$U$, which agree on $\gL$. Therefore, they agree everywhere.
\end{proof}

\section*{Appendix: Meyer sets in locally compact Abelian groups}

Let $G$ be a locally compact Abelian (LCA) group and let $\gL$ be a
Delone subset of $G$. We wish to compare the following two properties
that $\gL$ may have.

\begin{itemize}
\item[\textbf{M1}] $\gL - \gL \subset \gL + F$ for some finite set $F$;
\item[\textbf{M2}] $\gL -\gL$ is uniformly discrete.
\end{itemize}

Often these properties are taken as (equivalent) characterizations of
Meyer sets, which they are for groups of the form $\RR^d$. In fact, it
is easy to see that \textbf{M1} always implies \textbf{M2}. The reverse
implication for $\RR^d$ was proved by Lagarias \cite{Lagarias}. Here,
we prove it for all compactly generated LCA groups. The proof goes in
two steps.  First, we prove that the group generated by $\gL$ is
finitely generated.  In fact, this is equivalent to saying that $G$ is
compactly generated.  After that, we can basically follow Lagarias'
proof in the more general setting.

An apparently weaker concept than uniform discreteness is weak uniform
discreteness:

\begin{definition}
    $S \subset G$ is\/ \emph{weakly uniformly discrete} if for each
    compact subset $K$ of\/ $G$ and for all $a\in G$, $\mathrm{card} (S
    \cap (a+K))$ is bounded by a constant that depends only on $K$
    $($not on $a)$.
\end{definition}

Remarkably, as we shall see, for a Delone subset $\gL$ of a compactly
generated group $G$, the difference set $\varDelta:=\gL -\gL$ is
uniformly discrete if and only if it is weakly uniformly discrete.

\begin{prop} \label{M1toM2}
   If a Delone set $\gL$ satisfies {\rm\textbf{M1}}, it also satisfies 
   {\rm \textbf{M2}}.
\end{prop}

\begin{proof} 
Since $\gL$ is a Delone set, it is locally finite, i.e., $\gL\cap K$
is a finite set (or empty), for any compact set $K\subset G$. To
establish the uniform discreteness of $\varDelta$, it is sufficient to
show that $0$ is an isolated point of $\varDelta - \varDelta$.  Using
\textbf{M1} twice, one has
\[
   0 \; \in \; \varDelta - \varDelta \; \subset \;
   (\gL + F ) - (\gL + F) \; = \;
   \varDelta + (F - F) \; \subset \;
   \gL + F'
\]
where $F' = F + F - F$ is still a finite set. Consequently, $\gL + F'$
is locally finite, and $0$ must be an isolated point of it, hence also
of $\varDelta -\varDelta$. This gives \textbf{M2}.
\end{proof}

\begin{prop} \label{compactlyGenerated}
   Let\/ $G$ be an LCA group and let $\gL$ be relatively dense in\/ $G$.
   Suppose that $\langle \gL \rangle$ $($the subgroup of $G$ generated
   by $\gL)$ is finitely generated. Then, $G$ is compactly generated.
\end{prop}

\begin{proof}  
Let $F$ be a finite set that generates $\langle \gL \rangle$. As $\gL$
is relatively dense, there is a compact set $K\subset G$ with $G= \gL
+K$. Then, $F\cup K$ is compact and generates $G$.
\end{proof}

\begin{lemma} \label{projectionLemma}
   Let $G$ be an LCA group of the form $G' \times T$ where $T$ is a
   compact group. Then, the projection mapping $G \longrightarrow G'$
   defined by this splitting maps locally finite sets to locally
   finite sets.
\end{lemma}

\begin{proof} 
Suppose that $S \subset G$ is locally finite, but its projection $P' $
is not.  Then, there exists a compact set $K \subset G'$ with $P' \cap
K$ infinite and we have $S \cap (K \times T)$ infinite, too.  This is
a contradiction because $K \times T$ is compact.
\end{proof}

\begin{prop} \label{finitelyGenerated}
   Let $G$ be compactly generated. Let $\gL \subset G$ be relatively
   dense and suppose that $\varDelta=\gL - \gL$ is locally finite.
   Then, $\langle \gL \rangle$ is finitely generated.
\end{prop}

\begin{proof} 
By the structure theorem for compactly generated LCA groups
\cite[Thm.~9.8]{HW}, $G$ is isomorphic to $\RR^m \times \ZZ^n \times
T$, where $T$ is compact. We identify $G$ with this group and so can
also view it as a subgroup of $\RR^m \times \RR^n \times T =: G'
\times T$. We shall use $(')$ to indicate the projection map of $G'
\times T$ onto $G'$.  By Lemma~\ref{projectionLemma}, $\varDelta'
\subset G'$ is locally finite.  Also, $\gL$ is locally finite, due to
the corresponding property of $\varDelta$.

Select a compact set $C \subset G$ so that $\gL +C = G$. Since the
projection of $C$ into $G'$ is compact, we can find $R>0$ so that $\gL
+ (B_R \times T) \supset \RR^m \times \ZZ^n \times T$, where $B_R$ is
the open ball of radius $R$ around $0$ in $\RR^m \times \RR^n$.
Increasing $R$ if necessary, we may assume
\[ 
   \gL + (B_R \times T) \; = \; \RR^m \times \RR^n \times T .
\] 
Consider $F := (\gL\cup (\gL -\gL)) \cap (B_{2R} \times T)$, which is
finite.  It is plain that $F\subset\langle\gL\rangle$.  We show that
$\langle \gL\rangle = \langle F \rangle$. In fact, $\gL$ is contained
in the semigroup generated by $F$.

Let $\lambda \in \gL$. If $\lambda \in B_{2R} \times T$, then $\lambda
\in F$, so suppose $\lambda \notin B_{2R}\times T$. We need to get
closer to $0$, using a point of $F$. To this end, let $B_R(u')$ be
the open ball of radius $R$ around $u'$ in $G'$, where $u'$ is taken
to be the unique point which is at distance $R$ from $\lambda'$ on the
line segment $[0,\lambda']$ in $\RR^m \times \RR^n$ joining $0$ to
$\lambda'$. Thus, $\lambda' \in \overline{B_R(u')} \setminus B_R(u') =
\partial B_R (u')$. Let $u := (u',0) \in G'\times T$.  By our choice
of $R$, we can write $u = \mu_1 +b$, where $\mu_1\in\gL, b\in B_R
\times T$, so $\mu_1' = u' -b' \in B_R(u')$.

We have (i) $\mu_1 \in \gL$, (ii) $|\mu_1'| < |\lambda'|$, (iii)
$|\lambda' - \mu_1'| <2R$, where $|.|$ is the standard
Euclidean norm in $\RR^{m+n}$. Also, $\lambda - \mu_1 \in (B_{2R}
\times T) \cap \varDelta \subset F$ and so we have
\[
   \lambda = f_1 + \mu_1, \quad \mbox{with}\quad 
   f_1 \in F,\, \mu_1 \in \gL, \, |\mu_1'| <|\lambda'| \ts . 
\]
We now continue inductively, getting
\[
   \lambda \; = \; f_1 + \ldots + f_k + \mu_k
\]
where $f_1, \dots , f_k \in F$, $\mu_k \in \gL$ and $|\mu_k'| <
|\mu_{k-1}'| < \dots < |\lambda'|$, until $|\mu_k'| < 2R$. This must
happen for some $k$ since $\gL'$ is locally finite and $\gL' \cap
\overline{B_{\lambda'}(0)}$ is finite.  Then, $\mu_k \in F$ and we
have shown that $\lambda \in \langle F \rangle$.
\end{proof}

\begin{theorem} \label{lagariasTheorem}
   Let $G$ be a compactly generated LCA group.  Suppose that $\gL$ is
   a relatively dense subset of $G$ and that $\gL -\gL$ is weakly
   uniformly discrete.  Then, $\gL$ satisfies\/ {\rm \textbf{M1}}.
\end{theorem}

The proof of this result is really not different from that given in
\cite{Lagarias}. The only things to notice are that the full strength
of uniform discreteness of $\gL -\gL$ is not required and that we are
no longer confined to real spaces.

\begin{proof} 
We may assume that $0 \in \gL$, translating $\gL$ if necessary. Let $L
:= \langle \gL \rangle$ which is finitely generated by
Proposition~\ref{finitelyGenerated}: say $\langle \gL\rangle = \langle
e_1, \ldots , e_s \rangle$. For all $x \in L$, $x = \sum_{i=1}^s a_i
e_i$, $a_i \in \ZZ$, though not necessarily uniquely.

Define $\|x\| = \min \{ \sum |a_i| : x = \sum a_i e_i \}$. This
defines a norm on $L$ (where the proof of the triangle inequality
requires a short calculation).  For each $N \in \NN$, define $F(N) :=
\{x \in L : \|x\| \le N\}$, which is clearly a finite set.

Let $K \subset G$ be a symmetric compact neighbourhood of $0$ so that,
for every $u \in G$, we have $(u+K )\cap \gL \ne \varnothing$, and
also that $G$ is generated as a group by $K$. The goal is now to
show the existence of finitely many ``stepping stones'', forming
a set $F$, such that any difference $x-y$ of points in $\gL$ lies
in $\gL + F$.

To this end, let 
\begin{eqnarray*}
   M &:=& \max \{ {\textrm {card}} ((u + 2K) \cap (\gL -\gL)) : u \in
   G \} , \\ m&:=& \max \{ \|u\| : u \in (\gL -\gL) \cap ( K+K -K) \}.
\end{eqnarray*}
The former exists on the basis of the weak uniform discreteness of
$\gL -\gL$.

Let $x, y \in \gL$ and let $v: = y-x$.  Now, $x \in \ell K := K+
\ldots + K$ ($\ell$ summands) for some $\ell$, and we may write
\[
     x = k_1 + \ldots + k_\ell \, ,  
    \quad {\textrm{for some}} \;k_i \in K\ts .
\]

Let $x_0 := x, x_1 := x - k_1, x_2 := x - k_1 -k_2, \dots, x_\ell = 0$
and define the parallel sequence $y_i := x_i +v$, $0 \le i \le \ell$. 
Then, $y_i - y_{i+1} = x_i - x_{i+1} = k_{i+1} \in K$, for all
$0 \le i \le \ell-1$.

Choose $p_i, q_i \in \gL$ with $p_i - x_i , q_i - y_i \in K$ with the
special choices $p_0 = x, q_0 = y, p_\ell = 0 = x_\ell$. Note that
$q_0 = y = x+v = x_0 +v = y_0$, so in particular $q_0 - y_0 = 0 \in
K$.  Then, for each $i \in \{ 0, 1, \dots, \ell \}$, one has
$q_i - p_i -v = q_i -y_i +x_i - p_i \in 2K$.  Thus,
\[
    V \; := \; \{q_i - p_i : i = 0, \ldots , \ell\} 
    \subset \{(v+2K) \cap (\gL -\gL) \},
\]
so ${\textrm{card}} (V) \le M$.

Similarly, 
\[ 
   p_i - p_{i+1} \; = \; (p_i -x_i) +(x_i -x_{i+1}) +(x_{i+1} - p_{i+1})
   \; \subset \; (K+K-K) \cap (\gL - \gL)
\]
so $\|p_i - p_{i+1}\| \le m$.

In the same way, $\|q_i - q_{i+1}\| \le m$, so 
\[ 
   \|q_i - p_i - (q_{i+1} - p_{i+1})\| \; \le \; 2m \ts .
\]
Along with the bound on the cardinality of $V$, this gives rise to 
\[     \|u - u'\| \; \le \; 2mM    \]
for all $u,u' \in V$.

Now, $v = q^{}_0 -p^{}_0 = y-x \in V$ and $q^{}_\ell = q^{}_\ell -
p^{}_\ell \in V$, so $\|v - q^{}_\ell\| \le 2mM$. Since $v, q_\ell \in
L$, $v-q_\ell \in F(2mM)$ and
\[
   y-x \; = \; v \in q_\ell + F(2mM) 
   \; \subset \; \gL + F(2mM) \ts .
\]
Since $x,y \in \gL$ were arbitrary, $\gL -\gL \subset \gL + F(2mM)$. 
\end{proof}

\begin{coro} 
  Let $G$ be a compactly generated LCA group. Let $\gL\subset G$ be a
  Delone set. Then, $\gL - \gL$ is uniformly discrete if and only if
  it is weakly uniformly discrete.
\end{coro}
\begin{proof} 
Use Theorems~\ref{lagariasTheorem} and~\ref{M1toM2}. 
\end{proof}

\subsection*{Acknowledgements}
It is our pleasure to thank Jeong-Yup Lee, Bernd Sing and Nicolae
Strungaru for discussions and comments on the manuscript.  This work
was supported by the German Research Council (DFG), within
the CRC 701, and the Natural
Sciences and Engineering Council of Canada (NSERC).


\begin{thebibliography}{999}

\small

\bibitem{BHP}
M.~Baake, J.~Hermisson and P.~A.~B.~Pleasants,
\textit{The torus parametrization of quasiperiodic LI classes},
J.\ Phys.\ A:\ Math.\ Gen.\ \textbf{30} (1997) 3029--3056;
{\tt mp\_arc/02-168}.

\bibitem{BL} 
M.~Baake and D.~Lenz, 
\textit{Dynamical systems on
translation bounded measures:\ Pure point dynamical and diffraction
spectra}, Ergod.\ Th.\ \& Dynam.\ Syst.\ \textbf{24} (2004) 1867--93;
{\tt math.DS/0302061}.

\bibitem{BL2}
M.~Baake and D.~Lenz,
\textit{Deformation of Delone dynamical systems and pure point 
diffraction},
J.\ Fourier Anal.\ Appl.\ \textbf{11} (2005) 125--150;
{\tt math.DS/0404155}.

\bibitem{BLM} 
M.~Baake and D.~Lenz, 
\textit{Fourier modules and dual cut and project schemes}, 
in preparation.

\bibitem{BM}
M.~Baake and R.~V.~Moody,
\textit{Weighted Dirac combs with pure point diffraction},
J.\ reine angew.\ Math.\ (Crelle) \textbf{573} (2004) 61--94;
{\tt math.MG/0203030}.

\bibitem{Kola}
M.~Baake and B.~Sing,
\textit{Kolakoski-$(3,1)$ is a (deformed) model set},
Can.\ Math.\ Bull.\ \textbf{47} (2004) 168--190;
{\tt math.MG/0203025}.

\bibitem{BF}
C.~Berg and G.~Forst,
\textit{Potential Theory on Locally Compact Abelian Groups},
Springer, Berlin (1975).

\bibitem{BD} 
G.~Bernuau and M.~Duneau, 
\textit{Fourier analysis of deformed model sets}, 
in:~\textit{Directions in Mathematical
Quasicrystals}, eds.\ M.\ Baake and R.\ V.\ Moody, 
CRM Monograph Series, vol.\ 13, AMS, Providence, RI (2000), 
pp.\ 43--60.

\bibitem{Bou} 
N.~Bourbaki,
\textit{Elements of Mathematics:\ General Topology}, 
Chapters 1--4 and 5--10,
 reprint, Springer, Berlin (1989).

\bibitem{Cowley}
J.~M.~Cowley,
\textit{Diffraction Physics}, 3rd ed., 
North-Holland, Amsterdam (1995).

\bibitem{Dworkin}
S.~Dworkin, 
\textit{Spectral theory and $X\!$-ray diffraction},
J.\ Math.\ Phys.\ \textbf{34} (1993) 2965--2967.

\bibitem{GdeL}
J.~Gil de Lamadrid and L.~N.~Argabright, 
\textit{Almost Periodic Measures},
Memoirs of the AMS, vol.\ \textbf{65}, no.\ 428, AMS, 
Providence, RI (1990). 

\bibitem{HW}
E.~Hewitt and K.~A.~Ross,
\textit{Abstract Harmonic Analysis I},
2nd ed., Springer, New York (1979).

\bibitem{Hof}
A.~Hof,
\textit{On diffraction by aperiodic structures},
Commun.\ Math.\ Phys.\ \textbf{169} (1995) 25--43.

\bibitem{Hof3}
A.~Hof, 
\textit{Diffraction by aperiodic structures}, 
in:\ \textit{The Mathematics of Long-Range Aperiodic Order},
ed.\ R.\thinspace V.\ Moody,
NATO-ASI Series C 489, Kluwer, Dordrecht (1997), pp.\ 239--268.

\bibitem{K}
J.~Kwapisz,
\textit{Geometric coincidence conjecture and pure
discrete spectrum for unimodular tiling spaces},
talk given at the 2004 Banff meeting on
\emph{Aperiodic Order:\ Dynamical Systems, Combinatorics,
and Operators}.

\bibitem{Lagarias}
J.~C.~Lagarias,
\textit{Meyer's concept of quasicrystal and quasiregular sets},
Commun.\ Math.\ Phys.\ \textbf{179} (1996) 365--376.

\bibitem{Lee}
J.-Y.\ Lee, 
\textit{Substitution Delone sets with pure point spectrum are 
model sets}, preprint (2005).

\bibitem{LM}
J.-Y.\ Lee and R.~V.~Moody, 
\textit{A characterization of multi-colour model sets}, 
Annales Inst.\ H\'enri Poincar\'e \textbf{7} (2006), 125-143.

\bibitem{LMS-1}
J.-Y.\ Lee, R.~V.~Moody and B.~Solomyak,
\textit{Pure point dynamical and diffraction spectra},
Annales H.\ Poincar\'{e} \textbf{3} (2002) 1003--1018;
{\tt mp\_arc/02-39}.

\bibitem{LMS-2}
J.-Y.\ Lee, R.~V.~Moody and B.~Solomyak,
\textit{Consequences of pure point diffraction spectra for
multiset substitution systems},
Discr.\ Comput.\ Geom.\ \textbf{29} (2003) 525--560.

\bibitem{Len}
D.~Lenz, 
\textit{Continuty of eigenfunctions of uniquely ergodic 
dynamical systems and intensity of Bragg peaks}, 
preprint {\tt math-ph/0608026}.

\bibitem{Loomis}
L.~H.~Loomis, 
\textit{An Introduction to Abstract Harmonic Analysis},  
Van Nostrand, Princeton, NJ (1953).

\bibitem{Moody}
R.~V.~Moody,
\textit{Model sets:\ A survey}, in:\ 
\textit{From Quasicrystals to More Complex Systems}, 
eds.\ F.\ Axel, F.\ D\'enoyer and J.\ P.\ Gazeau,
EDP Sciences, Les Ulis, and
Springer, Berlin (2000), pp.\ 145--166;
{\tt math.MG/0002020}.

\bibitem{Moody2001}
R.~V.~Moody,
\textit{Uniform distribution in model sets},
Can.\ Math.\ Bulletin \textbf{45} (2002) 123--130.

\bibitem{MoodyLisbon}
R.~V.~Moody,
\textit{Mathematical quasicrystals:\ a tale of two topologies},
in:\ \textit{XIVth International Congress of Mathematical
Physics}, ed.\ J.-C.\ Zambrini,
World Scientific, Singapore (2005), pp.\ 68--77.

\bibitem{MoodyOberwolfach}
R.~V.~Moody,
\textit{The mathematics of aperiodic order},
Oberwolfach Reports \textbf{1} (2004) 1195--1198.

\bibitem{MS} 
R.~V.~Moody and N. Strungaru, 
\textit{Point sets and dynamical systems in the autocorrelation topology}, 
Can.\ Math.\ Bulletin \textbf{47} (2004) 82--99.

\bibitem{Ped}
G.~K.~Pedersen, 
\textit{Analysis Now}, 
Springer, New York (1989); rev.\ printing (1995).

\bibitem{Quef}
M.~Queff\'elec, 
\textit{Substitution Dynamical Systems -- Spectral Analysis},
Lecture Notes in Mathematics 1294, Springer, Berlin (1987).

\bibitem{Q}
B.~von Querenburg,
\textit{Mengentheoretische Topologie}, 2nd ed., 
Springer, Berlin (1979).

\bibitem{Rud}
 W.~Rudin,
\textit{Fourier Analysis on Groups}, 
 Wiley, New York (1962); reprint (1990).

\bibitem{Rob}
A.~E.~Robinson,  
\textit{On uniform convergence in the Wiener-Wintner theorem},  
J.\  London \ Math.\ Soc. \textbf{49} (1994) 493--501.

\bibitem{Shechtman}
D.~Shechtman, I.~Blech, D.~Gratias and J.\ts W.~Cahn,
\textit{Metallic phase with long-range orientational order
and no translation symmetry},
Phys.\ Rev.\ Lett.\ \textbf{53} (1984) 183--185.

\bibitem{Martin1}
M.~Schlottmann,
\textit{Cut-and-project sets in locally compact Abelian groups}, in:\
\textit{Quasicrystals and Discrete Geometry}, ed.\ J.\ Patera,
Fields Institute Monographs, vol.\ 10,
AMS, Providence, RI (1998), pp.\ 247--264.

\bibitem{Martin2}
M.~Schlottmann,
\textit{Generalized model sets and dynamical systems}, in:\
\textit{Directions in Mathematical Quasicrystals},
eds.\ M.\ Baake and R.\ V.\ Moody, CRM Monograph Series,
vol.\ 13, AMS, Providence, RI (2000), pp.\ 143--159.

\bibitem{SW}
B.~Sing and T.\ts R.~Welberry,
\textit{Deformed model sets and distorted Penrose tilings},
Z.\ Krist.\ \textbf{221} (2006) 621--634;
\texttt{mp\_arc/06-199}.

\bibitem{Boris}
B.~Solomyak, 
\textit{Spectrum of dynamical systems arising from Delone sets},
in:\ \textit{Quasicrystals and Discrete Geometry}, ed.\ J.\ Patera,
Fields Institute Monographs, vol.\ 10, 
AMS, Providence, RI (1998), pp.\ 265--275.

\bibitem{Boris2}
B.~Solomyak,
\textit{Dynamics of self-similar tilings},
Ergod.\ Th.\ \& Dynam.\ Syst.\ \textbf{17} (1997) 695--738;
Erratum: Ergod.\ Th.\ \& Dynam. \ Syst.\ \textbf{19} (1999) 1685. 
\bibitem{Wal} 
P.~Walters, 
\textit{An Introduction to Ergodic Theory},
Springer, New York (1982).

\bibitem{RS}
H.~Reiter and J.~D.~Stegeman,  
\textit{Classical Harmonic Analysis and Locally Compact Groups},
LMS Monographs, Clarendon Press, Oxford (2000).

\end{thebibliography}
\end{document}